\documentclass[11pt]{article}%
\usepackage{amsmath}
\usepackage{graphicx}%
\usepackage{amsfonts}%
\usepackage{amssymb}

\newcommand{\eqnb}{\begin{equation}}
\newcommand{\eqne}{\end{equation}}

\newtheorem{The}{Theorem}

\newtheorem{Def}{Definition}

\newtheorem{Lem}{Lemma}

\newtheorem{Rem}{Remark}

\usepackage{url}
\urldef{\mailsa}\path|liquanlin@tsinghua.edu.cn|
\urldef{\mailsb}\path|{qianzhiyong0926,fanruina}@stumail.ysu.edu.cn|

\begin{document}

\title{Fluid and Diffusion Limits for Bike Sharing Systems}
\author{Quan-Lin Li, Zhi-Yong Qian and Rui-Na Fan\\School of Economics and Management Sciences \\Yanshan University, Qinhuangdao 066004, P.R. China\\
\mailsa\\
\mailsb\\}
\maketitle

\begin{abstract}
Bike sharing systems have rapidly developed around the world, and
they are served as a promising strategy to improve urban traffic
congestion and to decrease polluting gas emissions. So far
performance analysis of bike sharing systems always exists many
difficulties and challenges under some more general factors. In this
paper, a more general large-scale bike sharing system is discussed
by means of heavy traffic approximation of multiclass closed
queueing networks with non-exponential factors. Based on this, the
fluid scaled equations and the diffusion scaled equations are
established by means of the numbers of bikes both at the stations
and on the roads, respectively. Furthermore, the scaling processes
for the numbers of bikes both at the stations and on the roads are
proved to converge in distribution to a semimartingale reflecting
Brownian motion (SRBM) in a $N^{2}$-dimensional box, and also the
fluid and diffusion limit theorems are obtained. Furthermore,
performance analysis of the bike sharing system is provided. Thus
the results and methodology of this paper provide new highlight in
the study of more general large-scale bike sharing systems.

\vskip                               0.5cm

\textbf{Keywords:} Bike sharing systems, fluid limit, diffusion
limit, semimartingale reflecting Brownian motion.

\end{abstract}

\section{Introduction}

Bike sharing systems have become an important way of urban
transportation due to its accessibility and affordability, and they
are widely deployed in more than 600 major cities around the world.
Bike sharing systems are regarded as promising solutions to reduce
congestion of traffic and parking, automobile exhaust pollution,
transportation noise, and so on. For some survey and development of
bike sharing systems, readers may refer to, DeMaio
\cite{DeMaio:2009}, Shaheen et al. \cite{Shaheen:2010}, Shu et al.
\cite{Shu:2013}, Labadi et al. \cite{Lab:2015}, and Meddin and
DeMaio \cite{Med:2012}.

Two major operational issues of bike sharing systems are to care for
(i) the non-empty: sufficient bikes parked at each station in order
to be able to rent a bike at any time; and (ii) the non-full:
suitable bike parking capacity designed for each station in order to
be able to return a bike in real time. Thus the empty or full
stations are called problematic stations. Up till now, efficient
measures are developed in the study of problematic stations,
including time-nonhomogeneous demand forecasting, average bike
inventory level, timely bike repositioning, and probability analysis
of problematic stations.

So far queueing models and Markov processes have been applied to
characterizing important steady-state performance of the bike
sharing systems. Important prior works on the bike sharing models
include the M/M/1/$C$ queue by Leurent \cite{Leurent:2012} and
Schuijbroek et al. \cite{Schuijbroek:2013}; the time-inhomogeneous
M$(t)$/M$(t)/1/C$ model by Raviv et al. \cite{Raviv:2013b} and Raviv
and Kolka \cite{Raviv:2013a}; the queueing networks by Kochel et al.
\cite{Kochel:2003}, Savin et al. \cite{Savin:2005}, Adelman
\cite{Adelman:2007}, George and Xia \cite{George:2010,George:2011}
and Li et al. \cite{Li:2016c}; the fluid models combining with
Markov decision processes by Waserhole and Jost
\cite{Waserhole:2012,Waserhole:2014}; the mean-field theory by
Fricker et al. \cite{Fricker:2012}, Fricker and Gast
\cite{Fricker:2014} and Fricker and Tibi \cite{Fricker:2015}; the
time-inhomogeneous M$(t)$/M$(t)/1/K$ and MAP$(t)$/MAP$(t)/1/K+2L+1$
queues combining with mean-field theory by Li et al. \cite{Li:2016a}
and Li and Fan \cite{Li:2016b}.

An important and realistic feature of bike sharing systems is the
time-varying arrivals of bike users and their random travel times.
In general, analysis of bike sharing systems with non-Poisson user
arrivals and general travel times are always very difficult and
challenging because more complicated multiclass closed queueing
networks are established to deal with bike sharing systems. See Li
et al. \cite{Li:2016c} for more interpretations. For this, fluid and
diffusion approximations may be an effective and better method in
the study of more general bike sharing systems. This motivates us in
this paper to develop fluid and diffusion limits for more general
large-scale bike sharing systems.

Fluid and diffusion approximations are usually applied to analysis
of more general large-scale complicated queueing networks, which
possibly originate in some practical systems including communication
networks, manufacturing systems, transportation networks and so
forth. See excellent monographs by, for example, Harrison
\cite{Harrison:1985}, Chen and Yao \cite{Chen:2001}, Whitt
\cite{Whitt:2002}. For the bike sharing system, further useful
information is introduced as follows. {\bf (a)} For heavy traffic
approximation of closed queueing networks, readers may refer to,
such as, Harrison et al. \cite{Harrison:1990} for a closed queueing
network with homogeneous customer population and infinite buffer.
Chen and Mandelbaum \cite{Chen:1991} for a closed Jackson network,
Harrison and Williams \cite{Harrison:1996} for a multiclass closed
network with two single-server stations and a fixed customer
population. Kumar \cite{Kumar:2000} for a two-server closed networks
in heavy traffic. {\bf (b)} For heavy traffic approximation of
queueing networks with finite buffers, important examples include,
Dai and Dai \cite{Dai:1999} obtained the SRBM of queue-length
process relying on a uniform oscillation result for solutions to a
family of Skorohod problems. Dai \cite{DaiW:1996} modeled the
queueing networks with finite buffers under a communication blocking
scheme, showed that the properly normalized queue length process
converges weakly to a reflected Brownian motion in a rectangular
box, and presented a general implementation via finite element
method to compute the stationary distribution of SRBM. Furthermore,
Dai \cite{Dai:2002} analyzed a multiclass queueing networks with
finite buffers and a feedforward routing structure under a blocking
scheme, and showed a pseudo-heavy-traffic limit theorem which stated
that the limit process of queue length is a reflecting Brownian
motion. {\bf (c)} There are some available results on heavy traffic
approximation of multiclass queueing networks, readers may refer to,
for instance, Harrison and Nguyen \cite{Harrison:1993}, Dai
\cite{Dai:1995}, Bramson \cite{Bramson:1998}, Meyn \cite{Meyn:2001}
and Majewski \cite{Majewski:2005}.

{\bf Contributions of this paper:} The main contributions of this
paper are threefold. The first contribution is to propose a more
general large-scale bike sharing system having renewal arrival
processes of bike users and general travel times, and to establish a
multiclass closed queueing network from the practical factors of the
bike-sharing system where bikes are abstracted as virtual customers,
while both stations and roads are regarded as virtual nodes or
servers. Note that the virtual customers (i.e. bikes) at stations
are of single class; while the virtual customers (i.e. bikes) on
roads are of two different classes due to two classes of different
bike travel or return times. The second contribution is to set up
the queue-length processes of the multiclass closed queueing network
through observing both some bikes parked at stations and the other
bikes ridded on roads. Such analysis gives the fluid scaled
equations and the diffusion scaled equations by means of the numbers
of bikes both at the stations and on the roads. The third
contribution is to prove that the scaling processes, corresponding
to the numbers of bikes both at the stations (having one class of
virtual customers) and on the roads (having two classes of virtual
customers), converge in distribution to a semimartingale reflecting
Brownian motion, and the fluid and diffusion limit theorems are
obtained in some simple versions. Based on this, performance
analysis of the bike sharing system is also given. Therefore, the
results and methodology given in this paper provide new highlight on
the study of more general large-scale bike sharing systems.

{\bf Organization of this paper:} The structure of this paper is
organized as follows. In Section 2, we describe a more general
large-scale bike sharing system with $N$ different stations and with
$N(N-1)$ different roads, while this system has renewal arrival
processes of bike users and general travel times on the roads. In
Section 3, we establish a multiclass closed queueing network from
practical factors of the bike-sharing system where bikes are
abstracted as virtual customers, while both stations and roads are
regarded as virtual nodes or servers. In Section 4, we set up the
queue-length processes of the multiclass closed queueing network by
means of the numbers of bikes both at the stations and on the roads,
and establish the fluid scaled equations and the diffusion scaled
equations. In Sections 5 and 6, we prove that the scaling processes
of the bike sharing system converge in distribution to a
semimartingale reflecting Brownian motion under heavy traffic
conditions, and obtain the fluid limit theorem and the diffusion
limit theorem, respectively. In Sections 7, we give performance
analysis of the bike sharing system by means of the fluid and
diffusion limits. Finally, some concluding remarks are described in
Section 8.

{\bf Useful notation:} We now introduce the notation used in the
paper. For positive integer $n$, the $n$-dimensional Euclidean space
is denoted by $\mathcal {R}^{n}$ and the $n$-dimensional positive
orthant is denoted by $\mathcal
{R}_{+}^{n}=\{x\in\mathcal{R}^{n}:x_{i}\geq0\}$. We definite
$D_{\mathcal{R}^{n}}[0,T]$ as the path space of all functions
$f:[0,T]\rightarrow \mathcal{R}^{n}$ which are right continuous and
have left limits. Define $\delta_{j,k}=1$ if $j=k$, else,
$\delta_{j,k}=0$. For a set $\mathcal {K}$, let $|\mathcal{K}|$
denote its cardinality. u.o.c. means that the convergence is
uniformly on compact set. A triple
$(\Omega,\mathcal{F},\{\mathcal{F}_{t}\})$ is called a filtered
space if $\Omega$ is a set, $\mathcal{F}$ is a $\sigma$-field of
subsets of $\Omega$, and $\{\mathcal{F}_{t},t\geq0\}$ is an
increasing family of sub-$\sigma $-fields of $\mathcal{F}$, i.e., a
filtration. If, in addition, $P$ is a probability measure on $\left(
\Omega,\mathcal{F}\right)  $, then
$(\Omega,\mathcal{F},\{\mathcal{F}_{t}\},P)$ is called a filtered
probability space. Let ${P_{x}}$ denote the unique family of
probability measures on $(\Omega,\mathcal{F})$, and $E_{x}$ be the
expectation operator under ${P_{x}}$.

\section{Model Description}\label{sec2:model}

In this section, we describe a more general large-scale bike sharing
system with $N$ different stations and with $N(N-1)$ different
roads, which has renewal arrival processes of bike users and general
travel times.

In the large-scale bike sharing system, a customer arrives at a
nonempty station, rents a bike, and uses it for a while, then he
returns the bike to a destination station and immediately leaves
this system. If a customer arrives at a empty station, then he
immediately leaves this system.

Now, we describe the bike sharing system including operations
mechanism, system parameters and mathematical notation as follows:

\textbf{(1) Stations and roads:} We assume that the bike sharing
system contains $N$ different stations and at most $N(N-1)$
different roads, where a pair of directed roads may be designed from
any station to another station. Also, we assume that at the initial
time $t=0$, each station has $C_{i}$ bikes and $K_{i}$ parking
positions, where $1\leq C_{i}\leq K_{i}<\infty$ for $i=1,\ldots,N$
and $\sum_{i=1}^{N}C_{i}>K_{j}$ for $j=1,\ldots,N$. Note that these
conditions make that some bikes can result in at least a full
station.

\textbf{(2) Arrival processes:} The arrivals of outside bike users
(or customers) at each station is a general renewal process. For
station $j$, let $u_{j}=\left\{ u_{j}(n),n\geq1\right\}$ be an
i.i.d. random sequence of exogenous interarrival times, where
$u_{j}(n)\geq0$ is the interarrival time between the $(n-1)$st
customer and the $n$th customer. We assume that $u_{j}(n)$ has the
mean $1/\lambda_{j}$ and the coefficient of variation $c_{a,j}$.

\textbf{(3) The bike return times:}

\textbf{(3.1) The first return:} Once an outside customer
successfully rents a bike from station $i$, then he rides on a road
directed to station $j$ with probability $p_{i\rightarrow j}$ for
$\sum_{j\neq i} ^{N}p_{i\rightarrow j}=1$, and his riding-bike time
$v_{i\rightarrow j}^{(1)}$ on the road $i\rightarrow j$ is a general
distribution with the mean $1/\mu_{i\rightarrow j}^{(1)}$ and the
coefficient of variation $c_{s,i\rightarrow j}^{(1)}$. If there is
at least one available parking position at station $j$, then the
customer directly returns his bike to station $j$, and immediately
leaves this system. Let $r^{i}=\{r_{j}^{i}(n),n\geq1\}$ be a
sequence of routing selections for $i,j=1,\ldots,N$ with $i\neq j$,
where $r_{j}^{i}(n)=1$ means that the $n$th customer rents a bike
from station $i$ and rides on a road directed to station $j$ (i.e.,
the customer rides on road $i\rightarrow j$), hence
$\Pr\{r_{j}^{i}(n)=1\}=p_{i\rightarrow j}$.

\textbf{(3.2) The second return:} From (3.1), if no parking position
is available at station $j$, then the customer has to ride the bike
to another station $l_{1}$ with probability $\alpha_{j\rightarrow
l_{1}}$ for $\sum_{l_{1}\neq j}^{N} \alpha_{j\rightarrow l_{1}}=1$,
and his riding-bike time $v_{j\rightarrow l_{1}}^{(2)}$ on road
$j\rightarrow l_{1}$ is also a general distribution with the mean
$1/\mu_{j\rightarrow l_{1}}^{(2)}$ and the coefficient of variation
$c_{s,j\rightarrow l_{1}}^{(2)}$. If there is at least one available
parking position at station $l_{1}$, then the customer directly
returns his bike and immediately leaves this bike sharing system.

\textbf{(3.3) The $(k+1)$st return for $k\geq2$:} From (3.2) and
more, we assume that this bike has not been returned at any station
yet through $k$ consecutive returns. In this case, the customer has
to try his $(k+1)$st lucky return, he will ride bike from the
$l_{k-1}$th full station to the $l_{k}$th station with probability
$\alpha_{l_{k-1}\rightarrow l_{k}}$ for $\sum_{l_{k}\neq l_{k-1}
}^{N}\alpha_{l_{k-1}\rightarrow l_{k}}=1$, and his riding-bike time
$v_{l_{k-1}\rightarrow l_{k}}^{(2)}$ on road $l_{k-1}\rightarrow
l_{k}$ is also a general distribution with the mean
$1/\mu_{l_{k-1}\rightarrow l_{k}}^{(2)}$ and the coefficient of
variation $c_{s,l_{k-1}\rightarrow l_{k}}^{(2)}$. If there is at
least one available parking station, then the customer directly
returns his bike and immediately leaves this bike sharing system;
otherwise he has to continuously try another station again. In the
next section, those bikes ridden under their first return are called
the first class of virtual customers; while those bikes ridden under
the $k$ ($k\geq2$) returns are called the second class of virtual
customers. Let $\bar{r}^{j}=\{\bar{r} _{i}^{j}(n),n\geq1\}$ be a
sequence of routing selections for $i,j=1,\ldots,N$ with $i\neq j$,
where $\bar{r}_{i}^{j}(n)=1$ means that the $n$th customer who can
not return the bike to the full station $j$ will deflect into road
$j\rightarrow i$, thus
$\Pr\{\bar{r}_{i}^{j}(n)=1\}=\alpha_{j\rightarrow i}$. Similarly,
let $r^{j\rightarrow i,(d)}=\{r^{j\rightarrow i,(d)}(n),n\geq1\}$ be
a sequence of routing selections for $i,j=1,\ldots,N$ with $i\neq
j$, $d=1,2$, where $r^{j\rightarrow i,(d)}(n)=1$ means the $n$th
customer of class $d$ who completes his short trip on road
$j\rightarrow i$ will return the bike to
station $i$, hence $\Pr\{r^{j\rightarrow i,(d)}%
(n)=1\}=p_{j\rightarrow i,i}=1$.

\textbf{(4) Two classes of riding-bike times:} In (3), there are two
classes of riding-bike times, who have two general distributions,
that is, there are two classes of virtual customers riding on each
road. Let $v_{j\rightarrow i}^{(d)}=\{v_{j\rightarrow i}
^{(d)}(n),n\geq1\}$ be a random sequence of riding-bike times of
class $d$ for $i,j=1,\ldots,N$ with $i\neq j$, $d=1,2$, where
$v_{j\rightarrow i}^{(d)}(n)$ is the riding-bike time for the $n$th
customer of class $d$ riding on the road $j\rightarrow i$. We assume
that $v_{j\rightarrow i}^{(d)}$ has the mean $1/\mu_{j\rightarrow
i}^{(d)}$ and the coefficient of variation $c_{s,j\rightarrow
i}^{(d)}$. To care for the expected riding-bike times, we set that
$\mu_{j\rightarrow i}^{(d)}=1/m_{j\rightarrow i}$ for $d=1$ and
$\mu_{j\rightarrow i}^{(d)}=1/\xi_{j\rightarrow i}$ for $d=2$.

\textbf{(5) The departure disciplines:} The customer departure has
two different cases: (a) an outside customer directly leaves the
bike sharing system if he arrives at an empty station; (b) if one
customer rents and uses a bike, and he finally returns the bike to a
station, then the customer completes his trip and immediately leaves
the bike sharing system.

For such a bike sharing system, Figure 1 outlines its physical
structure and associated operations.

\begin{figure}
\centering
\includegraphics[height=6.2cm]{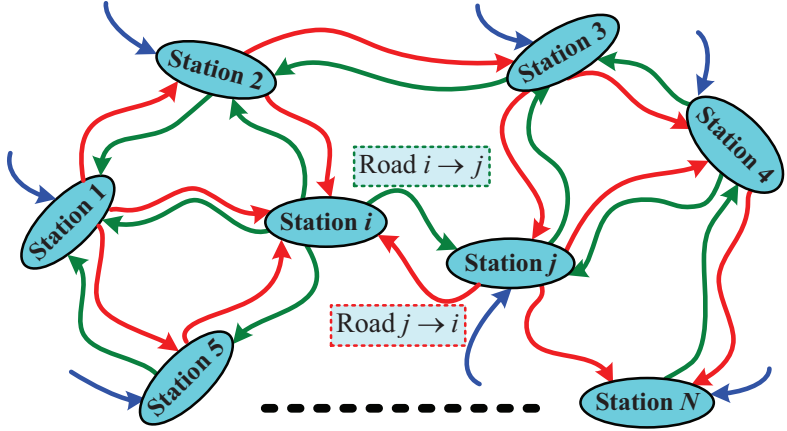}
\caption{The physical structure of the bike sharing system.}
\label{fig:phy-stru}
\end{figure}

\section{The Closed Queueing Network}\label{sec3:closed network}

In this section, we establish a multiclass closed queueing network
from the bike-sharing system where bikes are abstracted as virtual
customers, and both stations and roads are regarded as virtual nodes
or servers. Specifically, the stations contains only one class of
virtual customers; while the roads can contains two classes of
virtual customers.

In the bike sharing system, there are $N$ stations and $N(N-1)$
roads, and each bike can not leave this system, hence, the total
number of bikes in this system is fixed as $\sum_{i=1}^{N}C_{i}$.
Base on this, such a system can be regarded as a closed queueing
network with multiclass customers due to two types of different
travel or return times.

Let $S_{i}$ and $R_{i\rightarrow j}$ denote station $i$ and road
$i\rightarrow j$, respectively. Let SN denote the set of nodes
abstracted by the stations, and RN the set of nodes abstracted by
the roads. Clearly SN$=\{S_{i},i=1,\ldots,N\}$ and
RN$=\{R_{i\rightarrow j}: i,j=1,\ldots,N$ with $i\neq j\}$. Let
$n_{j}$ and $n_{i\rightarrow j}^{(d)}$ denote the numbers of bikes
parking in the $j$th station node and of bikes of class $d$ riding
on the road $i\rightarrow j$ node, respectively.

\textbf{(1) Virtual nodes:} Although the stations and the roads have
different physical attributes, they are all regarded as abstract
nodes in the closed queueing network.

\textbf{(2) Virtual customers:} The virtual customers are abstracted
by the bikes, which are either parked in the stations or ridden on
the roads. It is seen that only one class of virtual customers are
packed in the station nodes; while two classes of different virtual
customers are ridden on the road nodes due to their different return
times.

\textbf{(3) The routing matrix $P$:} To express the routing matrix,
we first define a mapping $\sigma(\cdot)$ as follow,
\[
\left\{
\begin{array}
[c]{l}%
\sigma(S_{i})=i\text{ \ \ \ \ \ \ \ \ \ \ \ for }i=1,\ldots,N,\\
\sigma(R_{i\rightarrow j})=i\langle j\rangle\text{ \ \ \ for
}i,j=1,\ldots,N,\text{ with }i\neq j.
\end{array}
\right.
\]
It is necessary to understand the mapping $\sigma(\cdot)$. For
example, $N=2$, $\sigma(S_{1})=1$,  $\sigma(S_{2})=2$,
$\sigma(R_{1\rightarrow 2})=1\langle 2 \rangle$,
$\sigma(R_{2\rightarrow 1})=2\langle 1 \rangle$, thus the routing
matrix is written as
\[
P=%
\begin{array}
[c]{cc} &
\begin{array}
[c]{cccc}%
1 & \text{ }2\text{ } & 1\langle 2\rangle & 2\langle 1\rangle
\end{array}
\\%
\begin{array}
[c]{c}%
1\\
2\\
1\langle 2\rangle\\
2\langle 1\rangle
\end{array}
& \left[
\begin{array}
[c]{cccc}%
\text{ \ \ \ } &  & \text{ \ \ \ } & \text{ \ \ \ }\\
& \text{ \ \ \ } &  & \\
&  &  & \\
&  &  &
\end{array}
\right].
\end{array}
\]
In this case, the component $p_{\tilde{i},\tilde{j}}$ of the routing
matrix $P$ denotes the probability
that a customer leaves node $\tilde{i}$ to node $\tilde{j}$, where%
\[
p_{\tilde{i},\tilde{j}}=\left\{
\begin{array}
[c]{l}%
1\text{ \ \ \ \ \ \ if }\tilde{i}=\sigma(R_{i\rightarrow j})\text{,
}\tilde
{j}=\sigma(S_{j})\text{,}\\
p_{i\rightarrow j}\text{ \ \ if }\tilde{i}=\sigma(S_{i})\text{,
}\tilde
{j}=\sigma(R_{i\rightarrow j})\text{,}\\
\alpha_{j\rightarrow k}\text{ \ if }\tilde{i}=\sigma(R_{i\rightarrow
j})\text{, }\tilde{j}=\sigma(R_{j\rightarrow k})\text{,}\\
0\text{ \ \ \ \ \ \ \ otherwise.}%
\end{array}
\right.
\]

\textbf{(4) The service processes in the station nodes:} For
$j\in\text{SN}$, the service process $S_{j}=\left\{
S_{j}(t),t\geq0\right\}$ of station node $j$, associated with the
interarrival time sequence $u_{j}=\left\{ u_{j}(n),n\geq1\right\}$
of the outside customers who arrive at station $j$, is given by
\[
S_{j}(t)=\sup\{n:U_{j}(n)\leq t\},
\]
where $U_{j}(n)=\sum_{l=1}^{n}u_{j}(l)$, $n\geq1$ and $U_{j}(0)=0$.
Let $b_{j}=\lambda_{j}1_{\{1\leq n_{j}\leq K_{j}\}}$.

\textbf{(5) The service processes in the road nodes:} For
$i,j=1,\ldots,N$ with $i\neq
j$ and $d=1,2$, the service process $S_{j\rightarrow i}^{(d)}%
=\{S_{j\rightarrow i}^{(d)}(t),t\geq0\}$ of road node $j\rightarrow
i$, associated with the riding-bike time sequence $v_{j\rightarrow
i}^{(d)}=\{v_{j\rightarrow i} ^{(d)}(n),n\geq1\}$ of the customers
of class $d$ ridden on road $j\rightarrow i$, is given by
\[
S_{j\rightarrow i}^{(d)}(t)=\sup\{n:V_{j\rightarrow i}^{(d)}(n)\leq
t\},
\]
where $V_{j\rightarrow i}^{(d)}(n)=\sum_{l=1}^{n}v_{j\rightarrow
i}^{(d)}(l),\text{ }n\geq1\text{ and }V_{j\rightarrow
i}^{(d)}(0)=0$. We write
\[
b_{j\rightarrow i}^{(d)}=n_{j\rightarrow i}^{(d)}\mu_{j\rightarrow i}%
^{(d)}=\left\{
\begin{array}
[c]{c}%
n_{j\rightarrow i}^{(1)}\frac{1}{m_{j\rightarrow i}}\text{ \ \ \ \ \ }d=1,\\
n_{j\rightarrow i}^{(2)}\frac{1}{\xi_{j\rightarrow i}}\text{ \ \ \ \
\ \ }d=2.
\end{array}
\right.
\]

\textbf{(6) The routing processes in the station nodes:}

\textit{Case one:} For $j\in\text{SN}$, the routing process
$R^{j}=\{R_{i}^{j},i\neq j,i=1,\ldots,N\}$ and
$R_{i}^{j}=\{R_{i}^{j}(n),n\geq1\}$, associated with the routing
selecting sequence $r^{i}=\{r_{j}^{i}(n),n\geq1\}$ of
station $j$, is given by%
\[
R^{j}(n)=\sum_{l=1}^{n}r^{j}(l)\text{ or }R_{i}^{j}(n)=\sum_{l=1}^{n}r_{i}^{j}(l),n\geq1,%
\]
and the $i$th component of $R^{j}(n)$ is $R_{i}^{j}(n)$ associated
with probability $p_{j\rightarrow i}$.

\textit{Case two:} For $j\in\text{SN}$, the routing process
$\bar{R}^{j}=\{\bar{R}_{i}^{j},i\neq j,i=1,\ldots,N\}$ and
$\bar{R}_{i}^{j}=\{\bar{R}_{i}^{j}(n),n\geq1\}$, associated with the
routing deflecting sequence $\bar{r}^{j}=\{\bar{r}
_{i}^{j}(n),n\geq1\}$ of station $j$, is given by
\[
\bar{R}^{j}(n)=\sum_{l=1}^{n}\bar{r}^{j}(l)\text{ or
}\bar{R}_{i}^{j}(n)=\sum_{l=1}^{n}\bar{r}_{i}^{j}(l),n\geq1,%
\]
and the $i$th component of $\bar{R}^{j}(n)$ is $\bar{R}_{i}^{j}(n)$
associated with probability $\alpha_{j\rightarrow i}$.

\textbf{(7) The routing processes in the road nodes:} For
$i,j=1,\ldots,N$ with $i\neq j$ and $d=1,2$, the routing process
$R^{j\rightarrow i,(d)}=\{R^{j\rightarrow i,(d)}(n),n\geq1\}$,
associated with the routing transferring sequence $r^{j\rightarrow
i,(d)}=\{r^{j\rightarrow i,(d)}(n),n\geq1\}$ of road $j\rightarrow
i$, is given by
\[
R^{j\rightarrow i,(d)}(n)=\sum_{l=1}^{n}r^{j\rightarrow
i,(d)}(l),n\geq1,%
\]
and the $R^{j\rightarrow i,(d)}(n)$ is associated with probability
$p_{j\rightarrow i,i}=1$.

\textbf{(8) Service disciplines:} The first come first served (FCFS)
discipline is assumed for all station nodes. A new processor sharing
(PS) is used for all the road nodes, where each customer of either
class one or class two is served by a general service time
distribution, as described in (4) and (5).

\section{The Joint Queueing Process}\label{sec4:fluid equations}

In this section, we set up the queue-length processes of the
multiclass closed queueing network by means of the numbers of bikes
both at the stations and on the roads, and establish the fluid
scaled equations and the diffusion scaled equations.

\textbf{(1) }$Q(t)=\{(Q_{j}(t),Q_{j\rightarrow i}^{(d)}(t)),i\neq
j,i,j=1,\ldots,N;d=1,2;t\geq0\}$, where $Q_{j}(t)$ and
$Q_{j\rightarrow i}^{(d)}(t)$ are the number of virtual customers at
station node $j$ and the numbers of virtual customers of class $d$
at the road $j\rightarrow i$ at time $t$, respectively.
Specifically, $Q_{j}(0)$ and $Q_{j\rightarrow i}^{(d)}(0)$ are the
number of virtual customers at station node $j$ and the number of
virtual customers of class $d$ on the road node $j\rightarrow i$ at
time $t=0$, respectively.

\textbf{(2) }$Y^{K}(t)=\{(Y_{j}^{K}(t)),j=1,\ldots,N;t\geq0\}$,
where $Y_{j}^{K}(t)$ is the cumulative number of virtual customers
deflecting from station node $j$ whose parking positions are full in
the time interval $\left[ 0,t\right]  $.

\textbf{(3) }$Y^{0}(t)=\{(Y_{j}^{0}(t),Y_{j\rightarrow
i}^{0,(d)}(t)),i\neq j,i,j=1,\ldots,N;d=1,2;t\geq0\}$, where
$Y_{j}^{0}(t)$ and $Y_{j\rightarrow i}^{0,(d)}(t)$ are the
cumulative amount of time that station node $j$ and the road node
$j\rightarrow i$ are idle (no available bike, i.e., empty) in the
time interval $[0,t]$, respectively.
\[
Y_{j}^{0}(t)=\int_{0}^{t}1\{Q_{j}(s)=0\}ds=t-B_{j}(t),
\]
\[
Y_{j\rightarrow i}^{0,(d)}(t)=\int_{0}^{t}1\{Q_{j\rightarrow i}^{(d)}%
(s)=0\}ds=t-B_{j\rightarrow i}^{(d)}(t).
\]

\textbf{(4) }$B(t)=\{(B_{j}(t),B_{j\rightarrow i}^{(d)}(t)),i\neq
j,i,j=1,\ldots,N;d=1,2;t\geq0\}$, where $B_{j}(t)$ and
$B_{j\rightarrow i}^{(d)}(t)$ are the cumulative amount of time that
the station node $j$ and the road node $j\rightarrow i$ are busy
(available bike, non-empty) in the time interval $[0,t]$,
respectively.
\[
B_{j}(t)=\int_{0}^{t}1\{0<Q_{j}(s)\leq K_{j}\}ds,
\]
\[
B_{j\rightarrow i}^{(d)}(t)=\int_{0}^{t}1\{Q_{j\rightarrow
i}^{(d)}(s)>0\}ds.
\]

\textbf{(5) }$B^{F}(t)=\{(B_{j}^{F}(t)),j=1,\ldots,N;t\geq0\}$,
where $B_{j}^{F}(t)$ is the cumulative amount of time that station
node $j$ is full (no available parking position) in the time
interval $[0,t]$,
\[
B_{j}^{F}(t)=\int_{0}^{t}1\{Q_{j}(s)=K_{j}\}ds.
\]

\textbf{(6) }$S_{j}(B_{j}(t))$ denotes the number of virtual
customers that have completed service at station node $j$ during the
time interval $[0,t]$; $S_{j\rightarrow i}^{(d)}(B_{j\rightarrow
i}^{(d)}(t))$ denotes the number of virtual customers of class $d$
that have completed service at road node $j\rightarrow i$ during the
time interval $[0,t]$.

\textbf{(7) }$R_{i}^{j} (S_{j}(B_{j}(t)))$ denotes the number of
virtual customers that enter station node $i$ (i.e., riding on road
$j\rightarrow i$) from station node $j$ during the time interval
$[0,t]$; $\bar{R}_{i}^{j}(Y_{j}^{K}(t))$ denotes the number of
virtual customers that enter station node $i$ from station $j$ whose
parking positions are full during the time interval $[0,t]$; and
$R^{j\rightarrow i,(d)}(S_{j\rightarrow i}^{(d)}(B_{j\rightarrow
i}^{(d)}(t)))$ denotes the number of virtual customers of class $d$
that enter station node $i$ from road node $j\rightarrow i$ during
the time interval $[0,t]$.

Now, we have the following flow balance relations for the station
nodes and the road nodes. For station node $j=1,\ldots,N$,
\begin{align}
Q_{j}(t) & =Q_{j}(0)+\sum_{d=1}^{2}\sum_{i\neq j}^{N}\left[
R^{i\rightarrow j,(d)}(S_{i\rightarrow j}^{(d)}(B_{i\rightarrow
j}^{(d)}(t)))-R^{i\rightarrow j,(d)}(S_{i\rightarrow
j}^{(d)}(B_{j}^{F}(t)))\right] \nonumber\\
& \text{\ \ \ \ } -S_{j}(B_{j}(t)). \label{eq1}
\end{align}
Note that $Y_{j}^{K}(t)=\sum_{d=1}^{2}\sum_{i\neq j}^{N}%
R^{i\rightarrow j,(d)}(S_{i\rightarrow j}^{(d)}(B_{j}^{F}(t)))$, we
have
\begin{equation}
Q_{j}(t) =Q_{j}(0)+\sum_{d=1}^{2}\sum_{i\neq j}^{N}R^{i\rightarrow
j,(d)}(S_{i\rightarrow j}^{(d)}(B_{i\rightarrow j}^{(d)}(t)))
-S_{j}(B_{j}(t))-Y_{j}^{K}(t). %
\label{eq2}
\end{equation}
For road node $j\rightarrow i$ for $i,j=1,\ldots,N$ with $i\neq j$
and $d=1,2$, we have
\begin{equation}
Q_{j\rightarrow i}^{(1)}(t)=Q_{j\rightarrow i}^{(1)}(0)+R_{i}^{j}(S_{j}%
(B_{j}(t)))-S_{j\rightarrow i}^{(1)}(B_{j\rightarrow i}^{(1)}(t)), \label{eq5}%
\end{equation}
\begin{equation}
Q_{j\rightarrow i}^{(2)}(t)=Q_{j\rightarrow i}^{(2)}(0)+\bar{R}_{i}^{j}%
(Y_{j}^{K}(t))-S_{j\rightarrow i}^{(2)}(B_{j\rightarrow
i}^{(2)}(t)). \label{eq6}
\end{equation}

Because the total number of bikes in this bike sharing system is
fixed as
$\sum_{i=1}^{N}C_{i}$, we get that for $t\geq0$%
\begin{equation}
\sum_{i=1}^{N}Q_{i}(t)+\sum_{d=1}^{2}\sum_{i\neq
j}^{N}Q_{i\rightarrow
j}^{(d)}(t)=\sum_{i=1}^{N}C_{i}. \label{eq7-2}%
\end{equation}
We now elaborate to apply a centering operation to the queue-length
representations of the station nodes and of the road nodes, and
rewrite (\ref{eq2}), (\ref{eq5}) and (\ref{eq6}) as follows:
\begin{equation}
Q(t)=X(t)+R^{0}Y^{0}(t)+R^{K}Y^{K}(t), \label{eq8}%
\end{equation}
where $X(t)=(X_{1}(t),X_{2}(t),\ldots,X_{N}(t))$, and $X_{j}(t)$ is given by%
\begin{align}
X_{j}(t)  &  =Q_{j}(0)+\sum_{d=1}^{2}\sum_{i\neq j}^{N}\left[
R^{i\rightarrow j,(d)}(S_{i\rightarrow j}^{(d)}(B_{i\rightarrow
j}^{(d)}(t)))-S_{i\rightarrow
j}^{(d)}(B_{i\rightarrow j}^{(d)}(t))\right] \nonumber\\
& \text{\ \ \ \ } +\sum_{d=1}^{2}\sum_{i\neq j}^{N}\left[  S_{i\rightarrow j}^{(d)}%
(B_{i\rightarrow j}^{(d)}(t))-b_{i\rightarrow j}^{(d)}B_{i\rightarrow j}%
^{(d)}(t)\right]  -\left[  S_{j}(B_{j}(t))-b_{j}B_{j}(t)\right] \nonumber\\
& \text{\ \ \ \ } -Y_{j}^{K}(t)+\theta_{j}t, \label{eq9-1}%
\end{align}
note that $R^{i\rightarrow j,(d)}(S_{i\rightarrow j}^{(d)}%
(B_{i\rightarrow j}^{(d)}(t)))=S_{i\rightarrow
j}^{(d)}(B_{i\rightarrow
j}^{(d)}(t))$, $X_{j}(t)$ is simplified as%
\begin{align}
X_{j}(t)  &  =Q_{j}(0)+\sum_{d=1}^{2}\sum_{i\neq j}^{N}\left[
S_{i\rightarrow j}^{(d)}(B_{i\rightarrow
j}^{(d)}(t))-b_{i\rightarrow j}^{(d)}B_{i\rightarrow
j}^{(d)}(t)\right] \nonumber\\
& \text{\ \ \ \ } -\left[  S_{j}(B_{j}(t))-b_{j}B_{j}(t)\right]
-Y_{j}^{K}(t)+\theta_{j}t,
\label{eq9}%
\end{align}%
\begin{equation}
\theta_{j}=\sum_{d=1}^{2}\sum_{i\neq j}^{N}b_{i\rightarrow
j}^{(d)}-b_{j},
\label{eq10}%
\end{equation}%
\begin{equation}
\left(  R^{0}Y^{0}(t)\right)_{\tilde{i},\tilde{j}}=\left\{
\begin{array}
[c]{l}%
b_{j}Y_{j}^{0}(t),\text{ \ \ \ \ \ \ \ \ \ \ \ \ \ \ \ \ \ \ \ if
}\tilde{i}=\sigma
(S_{i})\text{, and }\tilde{j}=\tilde{i},\\
-\sum_{d=1}^{2}b_{i\rightarrow j}^{(d)}Y_{i\rightarrow
j}^{0,(d)}(t),\text{ \ if }\tilde{i}=\sigma(S_{i})\text{, and
}\tilde{j}=\sigma(R_{i\rightarrow
j})\text{,}\\
0,\text{ \ \ \ \ \ \ \ \ \ \ \ \ \ \ \ \ \ \ \ \ \ \ \ \ \ \ \ \
otherwise,}
\end{array}
\right.  \label{eq11}%
\end{equation}%
\begin{equation}
(R^{K}Y^{K}(t))_{\tilde{i},\tilde{j}}=\left\{
\begin{array}
[c]{l}%
-Y_{j}^{K}(t),\text{ \ if }\tilde{i}=\sigma
(S_{i})\text{, and }\tilde{i}=\tilde{j},\\
0,\text{ \ \ \ \ \ \ \ \ \ \ otherwise. }%
\end{array}
\right.  \label{eq12}%
\end{equation}
For road node $j\rightarrow i$ $(i,j=1,\ldots,N$ with $i\neq j$ and
$d=1,2)$, $X_{j\rightarrow i}^{(d)}(t)$ is given by,%
\begin{align}
X_{j\rightarrow i}^{(1)}(t)  &  =Q_{j\rightarrow i}^{(1)}(0)+\left[  R_{i}%
^{j}(S_{j}(B_{j}(t)))-p_{j\rightarrow i}S_{j}(B_{j}(t))\right] \nonumber\\
& \text{\ \ \ \ } +\left[  p_{j\rightarrow
i}(S_{j}(B_{j}(t))-b_{j}B_{j}(t))\right]
\nonumber\\
& \text{\ \ \ \ } -\left[  S_{j\rightarrow i}^{(1)}(B_{j\rightarrow i}^{(1)}%
(t))-b_{j\rightarrow i}^{(1)}B_{j\rightarrow i}^{(1)}(t)\right]
+\theta_{j\rightarrow i}^{(1)}t, \label{eq13}%
\end{align}%
\begin{equation}
\theta_{j\rightarrow i}^{(1)}=p_{j\rightarrow i}b_{j}-b_{j\rightarrow i}%
^{(1)}, \label{eq14}%
\end{equation}%
\begin{equation}
\left(  R^{0}Y^{0}(t)\right)  _{\tilde{i},\tilde{j}}=\left\{
\begin{array}
[c]{l}%
b_{j\rightarrow i}^{(1)}Y_{j\rightarrow i}^{0,(1)}(t),\text{ \ if }%
\tilde{i}=\sigma(R_{j\rightarrow i})\text{ and }\tilde{j}=\tilde{i},\\
-p_{j\rightarrow i}b_{j}Y_{j}^{0}(t),\text{ if
}\tilde{i}=\sigma(R_{j\rightarrow i})\text{ and
}\tilde{j}=\sigma(S_{j}),\\
0,\text{ \ \ \ \ \ \ \ \ \ \ \ \ \ \ \ \ \ otherwise, }
\end{array}
\right.  \label{eq15}%
\end{equation}%
\begin{equation}
(R^{K}Y^{K}(t))_{\tilde{i},\tilde{j}}=0.\label{eq16}%
\end{equation}
\begin{align}
X_{j\rightarrow i}^{(2)}(t)  &  =Q_{j\rightarrow i}^{(2)}(0)+\left[
\bar {R}_{i}^{j}(Y_{j}^{K}(t))-\alpha_{j\rightarrow
i}Y_{j}^{K}(t)\right]
\nonumber\\
& \text{\ \ \ \ } -\left[  S_{j\rightarrow i}^{(2)}(B_{j\rightarrow i}^{(2)}%
(t))-b_{j\rightarrow i}^{(2)}B_{j\rightarrow i}^{(2)}(t)\right]
+\theta_{j\rightarrow i}^{(2)}t, \label{eq17}%
\end{align}%
\begin{equation}
\theta_{j\rightarrow i}^{(2)}=-b_{j\rightarrow i}^{(2)}, \label{eq18}%
\end{equation}%
\begin{equation}
\left(  R^{0}Y^{0}(t)\right)  _{\tilde{i},\tilde{j}}=\left\{
\begin{array}
[c]{l}%
b_{j\rightarrow i}^{(2)}Y_{j\rightarrow i}^{0,(2)}(t),\text{ \ if }%
\tilde{i}=\sigma(R_{j\rightarrow i})\text{ and }\tilde{j}=\tilde{i},\\
0,\text{ \ \ \ \ \ \ \ \ \ \ \ \ \ \ \ \ otherwise,}
\end{array}
\right.  \label{eq19}%
\end{equation}%
\begin{equation}
(R^{K}Y^{K}(t))_{\tilde{i},\tilde{j}}=\left\{
\begin{array}
[c]{l}%
\alpha_{j\rightarrow i}Y_{j}^{K}(t),\text{ \ if
}\tilde{i}=\sigma(R_{j\rightarrow
i})\text{ and }\tilde{j}=\sigma(S_{j}),\\
0,\text{ \ \ \ \ \ \ \ \ \ \ \ \ \ \ otherwise, }%
\end{array}
\right.  \label{eq20}%
\end{equation}
For $i,j=1,\ldots,N\text{ with }i\neq j,\text{ and }d=1,2$,
$Q_{j}(t)$, $Q_{j\rightarrow i}^{(d)}(t)$, $Y_{j}^{0}(t)$,
$Y_{j}^{K}(t)$, $Y_{j\rightarrow i}^{0,(d)}(t)$ have some important
properties as follows:
\begin{equation}
0\leq Q_{j}(t)\leq K_{j}\text{; \ }0\leq Q_{j\rightarrow
i}^{(d)}(t)\leq
\sum_{i=1}^{N}C_{i}\text{; \ }t\geq0,\label{eq24}%
\end{equation}%
\begin{equation}
Y_{j}^{0}(0)=0\text{, }Y_{j}^{0}(t)\text{ is continuous and
nondecreasing,
}\label{eq25}%
\end{equation}%
\begin{equation}
Y_{j}^{K}(0)=0\text{, }Y_{j}^{K}(t)\text{ is continuous and
nondecreasing,
}\label{eq26}%
\end{equation}%
\begin{equation}
Y_{j\rightarrow i}^{0,(d)}(0)=0\text{, }Y_{j\rightarrow
i}^{0,(d)}(t)\text{ is
continuous and nondecreasing, }\label{eq27}%
\end{equation}%
\begin{equation}
Y_{j}^{0}(t)\text{ increases at times }t\text{ only when
}Q_{j}(t)=0\text{,
}\label{eq28}%
\end{equation}%
\begin{equation}
Y_{j}^{K}(t)\text{ increases at times }t\text{ only when }Q_{j}(t)=K_{j}%
\text{, }\label{eq29}%
\end{equation}%
\begin{equation}
Y_{j\rightarrow i}^{0,(d)}(t)\text{ increases at times }t\text{
only when }Q_{j\rightarrow i}^{(d)}(t)=0\text{.}\label{eq30}%
\end{equation}

In the remainder of this section, we provide a lemma to prove that
the matrix $R=(R^{0},R^{K})$ is an $\mathcal{S}$ - matrix, which
plays a key role in discussing existence and uniqueness of the SRBM
through the box polyhedron for the closed queueing network. Note
that $R^{0}$ and $R^{K}$ are defined in (\ref{eq15}) and
(\ref{eq16}) for $d=1$, and in (\ref{eq19}) and (\ref{eq20}) for
$d=2$. Also, the $i$th column of $R$ is denoted as the vector
$v_{i}.$ To analyze the matrix $R$, readers may refer to Theorem 1.3
in Dai and Williams \cite{Dai:1996} for more details.

The following definition comes from Dai and Williams
\cite{Dai:1996}, here we restate it for convenience of readers.

\begin{Def}
\label{dft:matrix}A square matrix A is called an $\mathcal{S}$ -
matrix if there is a vector $x\geq0$ such that $Ax>0$. The matrix A
is completely - $\mathcal{S}$ if and only if each principal
submatrix of A is an $\mathcal{S}$ - matrix.
\end{Def}

Notice that the capacity of station nodes is finite and the total
number of bikes in this bike sharing system is a fixed constant.
Without loss of generality, we assume that the capacity of each road
node is also finite, and the maximal capacity of each road is
$\sum_{i=1}^{N}C_{i}$ due to the fact that the total number of bikes
in this bike sharing system is $\sum_{i=1}^{N}C_{i}$. Therefore, the
state space $S$ of this close queueing network is a
$N^{2}$-dimensional box space with $2N^{2}$ boundary faces $F_{i}$,
given by
\begin{equation}
S\equiv\{x=(x_{1},\ldots,x_{N^{2}})^{^{\prime}}\in\mathcal{R}_{+}^{N^{2}}:0\leq
x_{i}\leq\sum_{i=1}^{N}C_{i}\}. \label{eq31}%
\end{equation}%
We write
\begin{equation}
F_{i}\equiv\{x\in S:x_{i}=0\},F_{i+N^{2}}\equiv\{x\in
S:x_{i}=K_{i}\}\text{
for }i\in\text{SN,} \label{eq32}%
\end{equation}%
\begin{equation}
F_{j}\equiv\{x\in S:x_{j}=0\},F_{j+N^{2}}\equiv\{x\in S:x_{j}=\sum_{i=1}%
^{N}C_{i}\}\text{ for }j\in\text{RN.} \label{eq32-2}%
\end{equation}
Let $J\equiv\{1,2,\ldots,2N^{2}\}$ be the index set of the faces,
and for each $\emptyset\neq\mathcal{K}\subset J$, define
$F_{\mathcal{K}}=\cap_{i\in\mathcal{K}}F_{i}$. We indicate that the
set $\mathcal{K} \subset J$ is maximal if
$\mathcal{K}\neq\emptyset$, $F_{\mathcal{K}} \neq\emptyset$, and
$F_{\mathcal{K}}\neq F_{\tilde{\mathcal{K}}}$ for any
$\mathcal{K}\subset\mathcal{\tilde{\mathcal{K}}}$ such that
$\mathcal{K} \neq\mathcal{\tilde{\mathcal{K}}}$. Thus, we can obtain
that the maximal set $\mathcal{K}$ is precisely the set of indexes
of $N^{2}$ distinct faces meeting at any vertex of $S$. Let
$\mathcal{N}$ be a $2N^{2}\times N^{2}$ matrix whose $i$th row is
given by the unit normal of face $F_{i},$ which directs to the
interior of $S$. We obtain,
\[
\mathcal{N}=\left[
\begin{array}
[c]{c}%
\begin{array}
[c]{cccc}%
\text{ }1\text{ } & 0 & \text{ }\cdots\text{ } & 0\\
0 & \text{ }1\text{ } & \cdots & 0\\
\cdot & \cdot & \cdots & \cdot\\
0 & 0 & \cdots & \text{ }1\text{ }%
\end{array}
\\
---------\\%
\begin{array}
[c]{cccc}%
-1\text{ } & 0 & \text{ }\cdots\text{ } & 0\\
0 & -1\text{ } & \cdots & 0\\
\cdot & \cdot & \cdots & \cdot\\
0 & 0 & \cdots & -1\text{ }%
\end{array}
\end{array}
\right]  .
\]

The state space $S$ has $2^{N^{2}}$ vertexes due to its box space
and each vertex given by
$(\cap_{i\in\alpha}F_{i})\cap(\cap_{i\in\beta}F_{i+N^{2}})$ for a
unique index set $\alpha\subset\{1,\ldots,N^{2}\}$ with $\beta
=\{1,\ldots,N^{2}\}\backslash\alpha$. Before we provide a lemma to
prove the $(NR)_{\mathcal{K}}$ (exactly $|\mathcal{K}|$ distinct
faces contain $F_{\mathcal{K}}$) is a special $\mathcal{S}$-matrix,
we give a geometric interpretation for a
$|\mathcal{K}|\times|\mathcal{K}|$ $\mathcal{S}$-matrix
$(NR)_{\mathcal{K}}$. At the each vertex of the box, we should make
sure that there is a positive linear combination
$x_{i}v_{i}+x_{j}v_{j+N^2}$, $x_{i}>0$ for $i\in\alpha$ and
$x_{j}>0$ for $j\in\beta$ such that $x_{i}v_{i}+x_{j}v_{j+N^{2}}$
directs to the interior of the state space $S$.

Now, we provide a lemma to indicate the matrix $(NR)_{\mathcal{K}}$
is an $\mathcal{S}$-matrix.
\begin{Lem}
\label{lem:promatrix} The matrix $(NR)_{\mathcal{K}}$ is an
$\mathcal{S}$-matrix for each maximal $\mathcal{K}\subset J$.
\end{Lem}
\textbf{Proof:} It is easy to check that
\[
\mathcal{N}R=\left(
\begin{array}
[c]{c}%
\text{ }R^{0}\text{ \ \ \ \ }R^{K}\\
-R^{0}\text{ \ }-R^{K}%
\end{array}
\right).
\]
Because the state space of the closed queueing network is a
$N^{2}$-dimensional box space, it has $2N^{2}$ faces. Now, let us
make a classify of those vertexes in this box space as follows:

\textbf{Type-1}: the vertexes are given by $(\cap_{i\in
A_{S}}F_{i})\cap(\cap_{j\in A_{R}}F_{j})$;

\textbf{Type-2}: the vertexes are given by $(\cap_{i\in
A_{S}}F_{i})\cap(\cap_{k\in B_{R}}F_{k})$;

\textbf{Type-3}: the vertexes are given by $(\cap_{l\in
B_{S}}F_{l})\cap(\cap_{j\in A_{R}}F_{j})$;

\textbf{Type-4}: the vertexes are given by $(\cap_{l\in
B_{S}}F_{l})\cap(\cap_{k\in B_{R}}F_{k})$;

\textbf{Type-5}: the vertexes are given by $(\cap_{i\in
A_{S}}F_{i})\cap(\cap_{j\in A_{R}}F_{j})\cap(\cap_{k\in
B_{R}\backslash A_{R}}F_{k})$;

\textbf{Type-6}: the vertexes are given by $(\cap_{l\in
B_{S}}F_{l})\cap(\cap_{j\in A_{R}}F_{j})\cap(\cap_{k\in
B_{R}\backslash A_{R}}F_{k})$;

\textbf{Type-7}: the vertexes are given by $(\cap_{j\in
A_{R}}F_{j})\cap(\cap_{i\in A_{S}}F_{i})\cap(\cap_{l\in
B_{S}\backslash A_{S}}F_{l})$;

\textbf{Type-8}: the vertexes are given by $(\cap_{k\in
B_{R}}F_{k})\cap(\cap_{i\in A_{S}}F_{i})\cap(\cap_{l\in
B_{S}\backslash A_{S}}F_{l}$);

\textbf{Type-9}: the vertexes are given by $(\cap_{i\in
A_{S}}F_{i})\cap(\cap_{l\in B_{S}\backslash
A_{S}}F_{l})\cap(\cap_{j\in A_{R}}F_{j})\cap(\cap_{k\in
B_{R}\backslash A_{R}}F_{k})$;

\noindent where $A_{S}$ and $A_{R}$ denote the set of index of face
$F_{i}=\{x_{i}=0\}$ for $i\in\text{SN}$ and $F_{j}=\{x_{j}=0\}$ for
$j\in\text{RN}$, respectively; $B_{S}$ and $B_{R}$ denote the set of
index of face $F_{l}=\{x_{l}=K_{l}\}$ for $l\in\text{SN}$ and
$F_{k}=\{x_{k}=\sum_{i=1}^{N}C_{i}+1\}$ for $k\in\text{RN}$,
respectively. According to the model description in Section
\ref{sec2:model}, it is seen that the following two cases can not be
established:

\textbf{Case 1}: All the station nodes are saturated when $1\leq
C_{i}<K_{i}<\infty$, namely, the reflection direction vector $v_{i}$
on face $F_{i}(i\in B_{S})$ can not simultaneously exist in the box
state space $S$ due to $\sum_{i=1}^{N}K_{i}>\sum_{i=1}^{N}C_{i}$.
Therefore, at the vertexes of type-3, there must be a positive
linear combination $x_{i}v_{i}+x_{j}v_{j}>0$ to direct to the
interior of state space $S$, where $x_{i}\geq0$ for $i\in A_{R}$ and
$x_{j}\geq0$ for $j\in B_{S}$.

\textbf{Case 2}: Any road node is full, namely, the faces
$F_{i}(i\in B_{R})$ does not have the reflection direction vector
$v_{i}$ in the box state space $S$. In other word, the reflection
direction vector $v_{i}$ on face $F_{i}$ $(i\in B_{R})$ is zero
vector. Therefore, at the vertexes of type-2, type-4, type-5,
type-6, type-8 and type-9, there must be a positive linear
combination who directs to the interior of state space $S$.

Now, we should only prove that at these vertexes of type-1, type-7
and type-3, where $C_{i}=K_{i}$, there also is a positive linear
combination who directs to the interior of the state space $S$.

At the vertexes of type-1, we only should prove that the matrix
$R^{0}$ in the matrix $\mathcal{N}R$ is an $\mathcal{S}$-matrix for
$d=1,2$. It is clear that the matrix $R^{0}$ is an
$\mathcal{S}$-matrix due to the fact that all the diagonal elements
of $R^{0}$ are positive.

At the vertexes of type-7 and of type-3, for $C_{i}=K_{i}$ and
$d=1,2$,\ we can
rewrite the $(\mathcal{N}R)_{\mathcal{K}}$ as the following form:%
\[
M=(\mathcal{N}R)_{\mathcal{K}}=\left(
\begin{array}
[c]{c}%
M_{1}\text{ \ }M_{2}\\
M_{1}\text{ \ }M_{4}%
\end{array}
\right)  =\left(
\begin{array}
[c]{c}%
M_{1}\text{ \ }0\\
0\text{ \ }M_{4}%
\end{array}
\right)  +\left(
\begin{array}
[c]{c}%
0\text{ \ }M_{2}\\
M_{3}\text{ \ }0
\end{array}
\right)  .
\]
where $M_{1}$ is a submatrix of $R^{0}$, which contains $i$th row
(column) and $i$th column (row) of $R^{0}$ simultaneously with
$i\in\alpha\subset \{1,\ldots,N^{2}\}$. Because the $R^{0}$ is a
complete $\mathcal{S}$-matrix, $M_{1}$ is an $\mathcal{S}$-matrix.
$M_{4}$ is also a submatrix of $-R^{K}$, which also contains
$i+N^{2}$th row (column) and $i+N^{2}$th column (row) of $-R^{K}$
simultaneously with $i\in\beta =\{1,\ldots,N^{2}\}\backslash\alpha$.
At the same time, $M_{4}$ is a diagonal matrix whose diagonal
element is unit one, hence $M_{4}$ is also an $\mathcal{S}$-matrix.
$M_{2}$ is a submatrix of $R^{K}$ and $M_{3}$ is a submatrix of
$-R^{0}$. Because $M_{2}$ and $M_{3}$ do not contain any diagonal
elements of $R^{K}$ and $-R^{0}$, $M_{2}$ and $M_{3}$ are both
nonnegative matrices. Therefore, there must be a positive linear
combination who direct to the interior of the state space $S$ at the
vertexes of type-7 and type-3, for $C_{i}=K_{i}$ and $d=1,2$. This
completes the proof. \textbf{{\rule{0.08in}{0.08in}}}

\section{Fluid Limits}\label{sec5:fluid limit}

In this section, we provide a fluid limit theorem for the queueing
processes of the closed queueing network corresponding to the bike
sharing system.

It follows from the functional strong law of large numbers
(FSLLN) that as $t\rightarrow\infty$%
\begin{equation}
(\frac{1}{t}S_{j}(t),\frac{1}{t}S_{j\rightarrow
i}^{(d)}(t))\rightarrow(b_{j},b_{j\rightarrow i}^{(d)}),\text{ \ }d=1,2, \label{eq34}%
\end{equation}
and as $n\rightarrow\infty$%
\begin{equation}
(\frac{1}{n}R_{i}^{j}(n),\frac{1}{n}\bar{R}_{i}^{j}(n),\frac{1}{n}R^{i\rightarrow
j,(d)}(n))\rightarrow(p_{j\rightarrow i},\alpha_{j\rightarrow
i},1),\text{ \ }d=1,2.
\label{eq37}%
\end{equation}

We consider a sequence of closed queueing networks, indexed by
$n=1,2,\ldots$, as described in Section \ref{sec3:closed network}.
Let $(\Omega^{n},\mathcal{F}^{n},P^{n})$ be the probability space on
which the $n$th closed queueing network is defined for the bike
sharing system. All the processes and parameters associated with the
$n$th network are appended with a superscript $n$.

For the $n$th network, the renewal service processes of the station
nodes and of the road nodes are expressed by
$S_{j}^{n}=\{S_{j}^{n}(t),t\geq0\}$ and $S_{j\rightarrow
i}^{(d),n}=\{S_{j\rightarrow i}^{(d),n}(t),t\geq0\}$, respectively.
Let $b_{j}^{n}$ and $b_{j\rightarrow i}^{(d),n}$ be the long run
average service rates of $S_{i}^{n}(t)$ and $S_{j\rightarrow
i}^{(d),n}(t)$, respectively. The vectors of the $N$ station
capacities and of their initial bike numbers are denoted as
$K^{n}=(K_{1}^{n},\ldots,K_{N}^{n})^{^{\prime}}$ and
$C^{n}=(C_{1}^{n},\ldots,C_{N}^{n})^{^{\prime}}$, respectively,
where $1\leq C_{i}^{n}\leq K_{N}^{n}<\infty$. For simplicity of
description, we write $R^{j,n}\text{ as }R^{j}$,
$\bar{R}^{j,n}\text{ as }\bar{R}^{j}$ and $R^{j\rightarrow
i,(d),n}\text{ as }R^{j\rightarrow i,(d)}$ for all $n\geq1$, i.e.,
the routing processes of the station nodes and of the road nodes are
compressed the number $n$. We append a superscript $n$ to the
performance indexes such as $Y_{j}^{0,n}(t)$, $Y_{j\rightarrow
i}^{0,(d),n}(t)$, $B_{j}^{n}(t)$ and $B_{j\rightarrow i}^{n}(t)$,
and the interesting processes $Q^{n} =((Q_{j}^{n}(t),Q_{j\rightarrow
i}^{(d),n}(t))^{^{\prime}}$ and $Y_{j}^{K,n}(t)$.

\textbf{The heavy traffic conditions: }We assume that as $n\rightarrow\infty$%
\begin{equation}
(b_{j}^{n},b_{j\rightarrow i}^{(d),n},\sqrt{n}\theta_{j}^{n},\sqrt{n}%
\theta_{j\rightarrow
i}^{(d),n},\frac{1}{\sqrt{n}}C_{i}^{n},\frac{1}{\sqrt{n}}K_{i}^{n})\rightarrow
(b_{j},b_{j\rightarrow i}^{(d)},\theta_{j},\theta_{j\rightarrow i}^{(d)}%
,C_{i},K_{i}),\label{eq38}%
\end{equation}
where $\theta_{j}^{n}=\sum_{d=1}^{2}\sum_{j\neq
i}^{N}b_{j\rightarrow i}^{(d),n}-b_{j}^{n}$; $\theta_{j\rightarrow
i}^{(1),n}=p_{j\rightarrow i}b_{j}^{n}-b_{j\rightarrow i}^{(1),n}$
and $\theta_{j\rightarrow i}^{(2),n}=-b_{j\rightarrow i}^{(2),n}$.
At the same time, we assume that for $i,j=1,\ldots,N$ with $i\neq
j$, $d=1,2$, all these limits are finite.

For the initial queue lengths $Q_{j}^{n}(0)$ and $Q_{j\rightarrow
i}^{(d),n}(0)$, we assume that as $n\rightarrow\infty$%
\begin{equation}
\bar{Q}_{j}^{n}(0)\equiv\frac{1}{n}Q_{j}^{n}(0)\rightarrow0\text{
and }\bar
{Q}_{j\rightarrow i}^{(d),n}(0)\equiv\frac{1}{n}Q_{j\rightarrow i}%
^{(d),n}(0)\rightarrow0. \label{eq42-2}%
\end{equation}
It follows from the functional strong law of large numbers that for
$d=1,2,$ as $n\rightarrow\infty$
\begin{align}
(\bar{S}_{j}^{n}(t),\bar{S}_{j\rightarrow i}^{(d),n}(t),\bar{R}_{i}%
^{j,n}(t), & \bar{\bar{R}}_{i}^{j,n}(t),\bar{R}^{j\rightarrow i,(d),n}%
(t)) \nonumber \\
& \rightarrow(b_{j}t,b_{j\rightarrow i}^{(d)}t,p_{j\rightarrow i}%
t,\alpha_{j\rightarrow i}t,t),\text{ u.o.c.,} \label{eq42-3}%
\end{align}
where%
\[
\bar{S}_{j}^{n}(t)=\frac{1}{n}S_{j}^{n}(nt)\text{,
}\bar{S}_{j\rightarrow
i}^{(d),n}(t)=\frac{1}{n}S_{j}^{(d),n}(nt)\text{, }\bar{R}_{i}^{j,n}%
(t)=\frac{1}{n}R_{i}^{j}(\left\lfloor nt\right\rfloor )\text{,}%
\]%
\[
\bar{\bar{R}}_{i}^{j,n}(t)=\frac{1}{n}\bar{R}_{i}^{j}(\left\lfloor
nt\right\rfloor )\text{, }\bar{R}^{j\rightarrow i,(d),n}(t)=\frac{1}%
{n}R^{j\rightarrow i,(d)}(\left\lfloor nt\right\rfloor ),
\]
and $\left\lfloor x \right\rfloor$ is the maximal integer part of
the real number $x$.

We give a notation: for any process $W^{n}=\{W^{n}%
(t),t\geq0\}$, we define its centered processes
$\hat{W}^{n}=\{\hat{W}^{n}(t),t\geq0\}$ by%
\[
\hat{W}^{n}(nt)=W^{n}(nt)-w^{n}nt,
\]
where $w^{n}$ is the mean of the process $W^{n}$.

For the station nodes and road nodes, we write some centered
processes as
\begin{equation}
\hat{S}_{j}^{n}(nt)=S_{j}^{n}(nt)-b_{j}^{n}nt,\text{
}\hat{S}_{j\rightarrow i}^{(d),n}(t)=S_{j\rightarrow
i}^{(d),n}(nt)-b_{j\rightarrow i}^{(d),n}nt,
\label{eq43}%
\end{equation}%
\begin{equation}
\hat{R}_{i}^{j,n}(t)=R_{i}^{j,n}(\left\lfloor nt\right\rfloor
)-p_{j\rightarrow i}\left\lfloor nt\right\rfloor) \text{, }\hat{\bar{R}}%
_{i}^{j,n}(t)=\bar{R}_{i}^{j,n}(\left\lfloor nt\right\rfloor
)-\alpha
_{j\rightarrow i}\left\lfloor nt\right\rfloor). \label{eq44}%
\end{equation}

For convenience of readers, we restate a lemma for the oscillation
result of a sequence of $(S^{n},R^{n})$-regulation problems in
convex polyhedrons, which is a summary restatement of Lemma 4.3 of
Dai and Williams \cite{Dai:1996} and the Theorem 3.1 of Dai
\cite{DaiW:1996}, whose proof is omitted here and can easily be
referred to Dai and Williams \cite{Dai:1996} and Dai
\cite{DaiW:1996} for more details.

This lemma prevails due to the fact that the state space of the box
polyhedron of this bike sharing system belongs to a simple convex
polyhedrons as analyzed in the last of Section 4. For a function $f$
defined from $\left[ t_{1},t_{2}\right] \subset\left[
0,\infty\right]  $ into $\mathcal {R}^k$ for some $k\geq1$, let
\[
Osc(f,\left[  t_{1},t_{2}\right]  )=\sup_{t_{1}\leq s\leq t\leq
t_{2}}\left| f(t)-f(s)\right|  .
\]

\begin{Lem}
\label{lem:oscillation} For any $T>0$, given a sequence of $\{x^{n}\}_{n=1}%
^{\infty}\in D_{\mathcal {R}^{N^{2}}}\left[  0,T\right]  $ with the
initial values $x^{n}(0)\in S^{n}$. Let $(z^{n},y^{n})$ be an
$(S^{n},R^{n})$-regulation of
$x^{n}$ over $\left[  0,T\right]  $, where $(z^{n},y^{n})\in D_{\mathcal{R}^{N^{2}}%
}\left[  0,T\right]  \times D_{\mathcal{R}_{+}^{2N^{2}}}\left[
0,T\right] $. Assuming that all $S^{n}$ have the same shape, i.e.,
the only difference is the corresponding boundary size $K_{i}^{n}$.
Assuming that $\{K_{i}^{n}\}$ belongs to some bounded set, and the
jump sizes of $y^{n}$ are bounded by $\Gamma^{n}$ for each $n$. Then
if $(\mathcal{N}R)_{\mathcal{K}}$ is an $\mathcal{S}$ -
matrix and $R^{n}\rightarrow R$ as $n\rightarrow\infty$, we have%
\[
Osc(z^{n},[t_{1},t_{2}])\leq C\text{ }\max\{Osc(x^{n},[t_{1},t_{2}%
]),\Gamma^{n}\},
\]%
\[
Osc(y^{n},[t_{1},t_{2}])\leq C\text{ }\max\{Osc(x^{n},[t_{1},t_{2}%
]),\Gamma^{n}\},
\]
where $C$ depends only on $(\mathcal{N},R,\left|
{\mathcal{K}}\right| )$ for all $\mathcal{K}\subset\Xi$, where $\Xi$
denotes the collection of subsets of $J\equiv \{1,2,\ldots,2N^{2}\}$
consisting of all maximal sets in $J$ together with the empty set.
\end{Lem}

\begin{The}
\label{thm:fluid}(Fluid Limit Theorem) Under Assumptions
(\ref{eq38}) to (\ref{eq42-3}), as $n\rightarrow\infty$, we have
\[
\left(  \bar{B}_{j}^{n}(t),\bar{B}_{j\rightarrow i}^{(d),n}(t),\bar{Y}_{j}%
^{0,n}(t),\bar{Y}_{j\rightarrow i}^{0,(d),n}(t)\right)
\rightarrow\left( \bar{\tau}_{j}(t),\bar{\tau}_{j\rightarrow
i}^{(d)}(t),\bar{Y}_{j}^{0}(t),\bar
{Y}_{j\rightarrow i}^{0,(d)}(t)\right)  \text{ \ u.o.c,}%
\]
where $\bar{\tau}_{j}(t)\equiv et$, $\bar{\tau}_{j\rightarrow i}%
^{(d)}(t)\equiv et$, $\bar{Y}_{j}(t)\equiv0$ and
$\bar{Y}_{j\rightarrow
i}^{(d)}(t)\equiv0$; $\bar{Y}_{j}^{0,n}(t)=\frac{1}{n}Y_{j}%
^{0,n}(nt)$, $\bar{Y}_{j\rightarrow
i}^{0,(d),n}(t)=\frac{1}{n}Y_{j\rightarrow i}^{0,(d),n}(nt)$,
$\bar{B}_{j}^{n}(t)=\frac{1}{n}B_{j}^{n}(nt)$ and $\bar
{B}_{j\rightarrow i}^{(d),n}(t)=\frac{1}{n}B_{j\rightarrow
i}^{(d),n}(nt)$ for $i,j=1,\ldots,N$ with $i\neq j$, $d=1,2$.
\end{The}
\textbf{Proof:} Recall the queue length process :
$Q(t)=X(t)+R^{0}Y^{0}(t)+R^{K}Y^{K}(t)$, where $X(t)$ is given by
(\ref{eq9}), (\ref{eq13}) and (\ref{eq17}) in Section 4. It follows
from (\ref{eq2}) to (\ref{eq6}) that the scaling queueing processes
for the station nodes and the road nodes are given by
\[
\bar{Q}^{n}(t)=\bar{Q}^{n}(0)+\bar{X}^{n}(t)+R^{0,n}\bar{Y}^{0,n}%
(t)+R^{K,n}\bar{Y}^{K,n}(t),
\]
where $\bar{Q}^{n}(t)=\frac{1}{n}Q^{n}(nt)$, $\bar{Q}^{n}(t)=\{(\bar{Q}%
_{j}^{n}(t),\bar{Q}_{j\rightarrow i}^{(d),n}(t)),i\neq
j,i,j=1,\ldots,N;d=1,2;t\geq0\}$;
$\bar{X}^{n}(t)=\frac{1}{n}X^{n}(nt)$,
$\bar{X}^{n}(t)=\{(\bar{X}_{j}^{n}(t),$ $\ \bar{X}_{j\rightarrow i}%
^{(d),n}(t)),i\neq j,i,j=1,\ldots,N;d=1,2;t\geq0\}$; $\bar
{Y}^{0,n}(t)=\frac{1}{n}Y^{0,n}(nt)$, $\bar{Y}^{0,n}(t)=\{(\bar{Y}_{j}%
^{0,n}(t),\bar{Y}_{j\rightarrow i}^{0,(d),n}(t)),i$ $\neq
j,i,j=1,\ldots,N;d=1,2;t\geq0\}$;
$\bar{Y}^{K,n}(t)=\frac{1}{n}Y^{K,n}(nt)$,
$\bar{Y}^{K,n}(t)=\{(\bar{Y}_{j}^{K,n}(t)),$ $j=1,\ldots,N\}$. For
each $n$, $\bar{Q}^{n}(t)$, $\bar{Y}^{n}(t)$ and $\bar
{Y}^{K,n}(t)$\ satisfy the properties (\ref{eq24}) to (\ref{eq30})
with the
state space $S^{n}$, given by%
\begin{align*}
S^{n}  & \equiv\left\{  x=(x_{1},\ldots,x_{N^{2}})^{^{\prime}}\in R_{+}%
^{N^{2}}:x_{i}\leq\bar{K}_{i}^{n}=\frac{K_{i}^{n}}{n}\text{ for
}i\in
\text{SN};\right.  \\
& \left.  \text{ \ \ \ \ \ and
}x_{i}\leq\frac{\sum_{i=1}^{N}C_{i}^{n}}{n}+1\text{ for
}i\in\text{RN}\right\}  .
\end{align*}
For station node $j=1,\ldots,N$, by using (\ref{eq2}), (\ref{eq9}),
(\ref{eq43}) and (\ref{eq44}), we have%
\begin{equation}
\bar{X}_{j}^{n}(t)\equiv\frac{1}{n}Q_{j}^{n}(0)+\frac{1}{n}\sum_{d=1}^{2}\sum_{i\neq
j}^{N}\hat{S}_{i\rightarrow j}^{(d),n}(n\bar{B}_{i\rightarrow
j}^{(d),n}(t))-\frac{1}{n}\hat{S}_{j}^{n}(n\bar{B}_{j}^{n}(t))+\frac{1}{n}
\theta_{j}^{n}nt. \label{eq45}%
\end{equation}
For road node $j\rightarrow i\text{ }(i,j=1,\ldots,N\text{ with
}i\neq j)$, by using
(\ref{eq13}), (\ref{eq17}), (\ref{eq43}) and (\ref{eq44}), we have,%
\begin{align}
\bar{X}_{j\rightarrow i}^{(1),n}(t)  &
\equiv\frac{1}{n}Q_{j\rightarrow i}^{(1),n}(0)+\frac{1}{n}\hat
{R}_{i}^{j,n}(n\bar{S}_{j}^{n}(\bar{B}_{j}^{n}(t)))\nonumber\\
& \text{\ \ \ } +\frac{1}{n}p_{j\rightarrow i}\hat{S}_{j}^{n}(n\bar{B}_{j}^{n}%
(t))-\frac{1}{n}\hat{S}_{j\rightarrow
i}^{(1),n}(n\bar{B}_{j\rightarrow
i}^{(1),n}(t))+\frac{1}{n}\theta_{j\rightarrow i}^{(1),n}nt, \label{eq46-1}%
\end{align}
\begin{align}
\bar{X}_{j\rightarrow i}^{(2),n}(t)\equiv\frac{1}{n}Q_{j\rightarrow
i}^{(2),n}(0)+\frac{1}{n}
\hat{\bar{R}}_{i}^{j,n}(n\bar{Y}_{j}^{K,n}(t))-\frac{1}{n}\hat{S}_{j\rightarrow
i}^{(2),n}(n\bar{B}_{j\rightarrow
i}^{(2),n}(t))+\frac{1}{n}\theta_{j\rightarrow i}^{(2),n}nt. \label{eq47}%
\end{align}
Note that \bigskip$\bar{B}_{j\rightarrow i}^{(1),n}(t)\leq t$, $\bar
{B}_{j}^{n}(t)\leq t$,
$\bar{Y}_{j}^{K,n}(t)\leq\sum_{i=1}^{N}C_{i}^{n}-K_{j}^{n}$, by
using (\ref{eq38}) to (\ref{eq42-3}) and the Skorohod Representation
Theorem, as $n\rightarrow\infty$, we have%
\[
\bar{X}^{n}(t)=(\bar{X}_{j}^{n}(t),\bar{X}_{j\rightarrow i}^{(d),n}%
(t))\rightarrow0,\text{ \ \ \ \ u.o.c.}
\]
Since the state space $S^{n}$ of this bike sharing system are the
boxes of the same shape in the $N^{2}$-dimensional space,
$(\mathcal{N}R)_{\mathcal{K}}$ is an $S$-$matrix$ and
$R^{n}\rightarrow R$ as $n\rightarrow\infty$. Then by Lemma
\ref{lem:oscillation} we have
\[
Osc(\bar{Y}^{0,n},[s,t]\subseteq\lbrack0,T])\leq C\text{ }Osc(\bar{X}%
^{n},[s,t]\subseteq\lbrack0,T]),
\]
for any $T\geq0$, where $C$ depends only on $R$ and $\mathcal{N}$
for $n$ large enough.%
\begin{align*}
0  &  \leq\lim_{n\rightarrow\infty}\inf
Osc(\bar{Y}^{0,n},[s,t]\subseteq
\lbrack0,T])\\
&  \leq\lim_{n\rightarrow\infty}\sup
Osc(\bar{Y}^{0,n},[s,t]\subseteq
\lbrack0,T])\\
&  \leq C\text{
}\lim_{n\rightarrow\infty}Osc(\bar{X}^{n},[s,t]\subseteq
\lbrack0,T])\\
&  =0,\text{ a.s.}
\end{align*}
where
$\bar{Y}^{n}(t)=(\bar{Y}_{j}^{0,n}(t)^{^{\prime}},\bar{Y}_{j\rightarrow
i}^{(d),0,n}(t)^{^{\prime}})^{^{\prime}}$. Notice that $Y^{n}(0)=0$
for all
$n$, we have%
\begin{equation}
\lim_{n\rightarrow\infty}\bar{Y}^{n}(t)=0,\text{ u.o.c.}
\label{eq50}%
\end{equation}
Since $\bar{B}_{j}^{n}(t)=t-\bar{Y}_{j}^{0,n}(t)$ and $\bar
{B}_{j\rightarrow i}^{(d),n}(t)=t-\bar{Y}_{j\rightarrow
i}^{(d),0,n}(t)$, we obtain the convergence of $\bar{B}_{j}^{n}(t)$
and $\bar{B}_{j\rightarrow i}^{(d),n}(t)$ for $i,j=1,\ldots,N$ with
$i\neq j$, $d=1,2$. This competes the proof.
\textbf{{\rule{0.08in}{0.08in}}}

\section{Diffusion limits}\label{sec6:diffusion limit}

In this section, we set up the diffusion scaled processes of the
queueing processes, and give their weak convergence results for the
multiclass closed queueing network corresponding to the bike sharing
system.

We introduce the diffusion scaling process for the process
$\hat{W}^{n}=\{\hat{W}^{n}(nt),t\geq0\}$, given by
\[
\tilde{W}^{n}(t)\equiv\frac{1}{\sqrt{n}}\hat{W}^{n}(nt)=\frac{1}{\sqrt{n}%
}(W^{n}(nt)-w^{n}nt).
\]

For the station nodes and the road nodes, we write
\begin{equation}
\tilde{S}_{j}^{n}(t)=\sqrt{n}\left(\frac{S_{j}^{n}(nt)}{n}-b_{j}^{n}t\right),\text{
 }
\tilde{S}_{j\rightarrow
i}^{(d),n}(t)=\sqrt{n}\left(\frac{S_{j\rightarrow i}
^{(d),n}(nt)}{n}-b_{j\rightarrow i}^{(d),n}t\right), \label{eq51-1}%
\end{equation}%
\begin{equation}
\tilde{R}_{i}^{j,n}(t)=\sqrt{n}\left(\frac{R_{i}^{j,n}(nt)}{n}-p_{j\rightarrow
i}t\right), \text{
 }\tilde{\bar{R}}_{i}^{j,n}(t)=\sqrt{n}\left(\frac{\bar{R}_{i}^{j,n}(nt)}{n}-\alpha_{j\rightarrow
i}t\right), \label{eq52-1}%
\end{equation}%
\begin{equation}
\tilde{R}^{j\rightarrow i,(d),n}(t)=\sqrt{n}\left(\frac{R^{j\rightarrow i,(d),n}(nt)}%
{n}-t\right). \label{eq52-3}%
\end{equation}

For the initial queueing processes $Q_{j}^{n}(0)$ and
$Q_{j\rightarrow i}^{n,(d)}(0)$ for $i,j=1,\ldots,N$ with $i\neq j$,
$d=1,2$, we assume that as $n\rightarrow
\infty$%
\begin{equation}
\tilde{Q}_{j}^{n}(0)\equiv\frac{1}{\sqrt{n}}Q_{j}^{n}(0)\Rightarrow\tilde
{Q}(0),
\end{equation}
\begin{equation}\tilde{Q}_{j\rightarrow i}^{(d),n}(0)\equiv\frac{1}{\sqrt
{n}}Q_{j\rightarrow i}^{(d),n}(0)\Rightarrow\tilde{Q}_{j\rightarrow i}%
^{(d)}(0). \label{eq53}%
\end{equation}

It follows from the Skorohod Representation Theorem and the
Donsker's Theorem that
\begin{align}
(\tilde{S}_{j}^{n}(t),\tilde{S}_{j\rightarrow i}^{(d),n}(t),\tilde
{R}_{i}^{j,n}(t), & \tilde{\bar{R}}_{i}^{j,n}(t),
\tilde{R}^{j\rightarrow i,(d),n}(t)) \nonumber \\
& \Rightarrow(\tilde{S}_{j}(t),\tilde{S}_{j\rightarrow i}%
^{(d)}(t),\tilde{R}_{i}^{j}(t),\tilde{\bar{R}}_{i}^{j}(t),\tilde
{R}^{j\rightarrow i,(d)}(t)),\label{eq53-1}%
\end{align}
where $\Rightarrow$ denotes weak convergence, and
$\tilde{S}_{j}(t)$, $\tilde{S}_{j\rightarrow i}^{(d)}(t)$, $\tilde
{R}_{i}^{j}(t)$, $\tilde{\bar{R}}_{i}^{j}(t)$ and
$\tilde{R}^{j\rightarrow i,(d)}(t)$ are all the Brownian motions
with drift zero and covariance matrices $\Gamma^{S}$,
$\Gamma^{R,S,l}$, $\Gamma^{\bar{R},S,l}$ and
$\Gamma^{R,S,j\rightarrow i}$, which are given by

\textbf{(1) }The covariance matrix of
$\tilde{S}(t)=(\tilde{S}_{j}(t),\tilde{S} _{j\rightarrow
i}^{(d)}(t))$ for $i,j=1,\ldots,N$ with $i\neq j$, $d=1,2$, is given
by
\[
\Gamma^{S}=\left(
\begin{array}
[c]{cc}%
\left(  \Gamma^{S,S}\right)  _{N\times N} & 0\\
0 & \left(  \Gamma^{S,R,(d)}\right)  _{\left(  N^{2}-N\right)
\times\left(
N^{2}-N\right)  }%
\end{array}
\right)  _{N^{2}\times N^{2}},
\]
where%
\[
\left(  \Gamma^{S,S}\right)  _{\tilde{i},\tilde{j}}=\left\{
\begin{array}
[c]{c}%
b_{i}c_{a,i}^{2}\delta_{\tilde{i},\tilde{j}},\text{ \ }\sigma(S_{i}%
)=\tilde{i}\text{,}\\
0,\text{ \ \ \ \ \ \ \ \ \ \ otherwise,}%
\end{array}
\right.
\]
\[
\left(  \Gamma^{S,R,(d)}\right)  _{\tilde{i},\tilde{j}}=\left\{
\begin{array}
[c]{l}%
b_{i\rightarrow j}^{(d)}(c_{s,i\rightarrow
j}^{(d)})^{2}\delta_{\tilde{\imath
},\tilde{j}},\text{ \ }\sigma(R_{i\rightarrow j})=\tilde{i}\text{,}\\
0,\text{ \ \ \ \ \ \ \ \ \ \ \ \ \ \ \ \ \ \ \ \ otherwise.}%
\end{array}
\right.
\]

\textbf{(2) }The covariance matrix of
$\tilde{R}(t)=(\tilde{R}^{l}(t))$  for $l=1,\ldots,N$, is given by
\[
\Gamma^{R,S,l}=\left(
\begin{array}
[c]{cc}%
0 & 0\\
0 & \left(  \Gamma^{R,S,l}\right)  _{\left(  N-1\right) \times\left(
N-1\right)  }%
\end{array}
\right)  _{N^{2}\times N^{2}},
\]
where%
\[
\left(  \Gamma^{R,S,l}\right)  _{\tilde{i},\tilde{j}}=\left\{
\begin{array}
[c]{l}%
p_{l\rightarrow k_{1}}(\delta_{\tilde{i},\tilde{j}}-p_{l\rightarrow
k_{2}}),\text{\ }\sigma(R_{l\rightarrow k_{1}})=\tilde{i},\text{
\ }\sigma(R_{l\rightarrow k_{2}})=\tilde{j}\text{,}\\
0,\text{ \ \ \ \ \ \ \ \ \ \ \ \ \ \ \ \ \ \ \ \ \ \ \ otherwise.}%
\end{array}
\right.
\]

\textbf{(3) }The covariance matrix of
$\tilde{\bar{R}}(t)=(\tilde{\bar{R}}^{l}(t))$ for $l=1,\ldots,N$, is
given by
\[
\Gamma^{\bar{R},S,l}=\left(
\begin{array}
[c]{cc}%
0 & 0\\
0 & \left(  \Gamma^{\bar{R},S,l}\right)  _{\left(  N-1\right)
\times\left(
N-1\right)  }%
\end{array}
\right)  _{N^{2}\times N^{2}},
\]
where%
\[
\left(  \Gamma^{\bar{R},S,l}\right) _{\tilde{i},\tilde{j}}=\left\{
\begin{array}
[c]{l}%
\alpha_{l\rightarrow k_{1}}(\delta_{\tilde{i},\tilde{j}}-\alpha
_{l\rightarrow k_{2}}),\text{ \ \ }\sigma(R_{l\rightarrow k_{1}}%
)=\tilde{i}\text{, }\sigma(R_{l\rightarrow k_{2}})=\tilde{j}\text{,}\\
0,\text{ \ \ \ \ \ \ \ \ \ \ \ \ \ \ \ \ \ \ \ \ \ \ \ \ \ otherwise.}%
\end{array}
\right.
\]

\textbf{(4) }The covariance matrix of
$\tilde{R}(t)=(\tilde{R}^{j\rightarrow
i,(d)}(t))$ for $i,j=1,\ldots,N$ with $i\neq j$, $d=1,2,$ is given by%
\[
\Gamma^{R,S,j\rightarrow i}=\left(
\begin{array}
[c]{cc}%
\left(  \Gamma^{R,R,j\rightarrow i}\right)  _{N\times N} & 0\\
0 & 0
\end{array}
\right)  _{N^{2}\times N^{2}},
\]
where%
\[
\left(  \Gamma^{R,R,j\rightarrow i}\right)
_{\tilde{l},\tilde{k}}=\left\{
\begin{array}
[c]{l}%
p_{j\rightarrow i,l}(\delta_{\tilde{l},\tilde{k}}-p_{j\rightarrow
i,k})=0,\text{ \ }\sigma(S_{l})=\tilde{l},\text{
}\sigma(S_{k})=\tilde
{k}\text{,}\\
0,\text{ \ \ \ \ \ \ \ \ \ \ \ \ \ \ \ \ \ \ \ \ \ \ \ \ \ \ \ \ \ \
otherwise.},
\end{array}
\right.
\]

Now, we prove adaptedness properties of the diffusion scaling
processes $(\tilde{Q}^{n}(t),$ $\tilde{X}^{n}(t),\tilde{Y}^{n}(t))$,
where $\tilde{Q}^{n}(t)=\frac{1}{\sqrt{n}}Q^{n}(nt)$, $\tilde
{Q}^{n}(t)=(\tilde{Q}_{j}^{n}(t),\tilde{Q}_{j\rightarrow
i}^{(d),n}(t))$; $\tilde{X}^{n}(t)=\frac{1}{\sqrt{n}}X^{n}(nt)$,
$\tilde{X}^{n}(t)=(\tilde
{X}_{j}^{n}(t),\tilde{X}_{j\rightarrow i}^{(d),n}(t))$; $\tilde{Y}%
^{0,n}(t)=\frac{1}{\sqrt{n}}Y^{0,n}(nt)$, $\tilde{Y}^{0,n}(t)=(\tilde{Y}%
_{j}^{0,n}(t)^{^{\prime}}$, $\tilde{Y}_{j\rightarrow i}^{0,(d),n}%
(t)^{^{\prime}})$,
$\tilde{Y}_{j}^{K,n}(t)=\frac{1}{\sqrt{n}}Y_{j}^{K,n}(nt)$.

Define
\begin{equation}
\varsigma_{t}^{n}=\sigma\{\tilde{Q}^{n}(0),\tilde{S}^{n}(s),\tilde{Y}%
^{0,n}(s),\tilde{Y}^{K,n}(t),s\leq t\},\label{eq54}%
\end{equation}
where $\tilde{Q}^{n}(0),\tilde{S}^{(d),n}(s)$,
$\tilde{R}^{(d),n}(s)$ and $\tilde{\bar{R}}^{n}(s)$\ are defined in
(\ref{eq51-1}) to (\ref{eq53}). Define
$T_{k}^{n}=(T_{k}^{j,n},T_{k}^{j\rightarrow i,(d),n})$, where
$T_{k}^{j,n}$ and $T_{k}^{j\rightarrow i,(d),n}$ denote the partial
sum of the service time sequence at station node $j$ and road node
$j\rightarrow i$,
respectively, for the $n$th network, that is,%
\[
T_{k}^{j,n}=\sum_{l=1}^{k}u_{j}^{n}(l),\text{  }T_{k}^{j\rightarrow
i,(d),n}=\sum_{l=1}^{k}v_{j\rightarrow i}^{(d),n}(l),
\]
with the initial condition $T_{0}^{n}\equiv0$. Notice that
$T_{k}^{n}=(T_{k}^{j,n},T_{k} ^{j\rightarrow i,(d),n})$ is a
$\varsigma_{t}^{n}-stopping$ $time$, and,
$0=T_{0}^{n}<T_{1}^{n}<T_{2}^{n}<\cdots<T_{k}^{n}\rightarrow\infty$
a.s. as $k\rightarrow\infty$ for each $n$ and $i,j=1,\ldots,N$ with
$i\neq j$, $d=1,2$. Let $\varsigma_{T_{k}^{(n)-}}$ denote the strict
past at time $T_{k}^{n}$. Then
\[
\varsigma_{T_{k}^{(n)-}}=\sigma(A_{t}\cap\{t<T_{k}^{n}\},A_{t}\in\varsigma
_{t}^{n},t\geq0).
\]
Because $T_{k}^{n}$ is a $\varsigma_{t}^{n}$-stopping time,
$u_{j}^{n}(k+1)$ and $v_{j\rightarrow i}^{(d),n}(k+1)$ are
independent of the history of the network before the time at which
the $k$th customer is served at station node $j$ and road node
$j\rightarrow i$. Therefore, $T_{k}^{n}$ is $\varsigma
_{T_{k}^{(n)-}}$-measurable, $u_{j}^{n}(k+1)$ is independent of
$\varsigma _{T_{k}^{(j,n)-}\text{ }}$, and $v_{j\rightarrow
i}^{(d),n}(k+1)$ is independent of $\varsigma_{T_{k}^{(j\rightarrow
i,(d),n)-}}$.
\begin{The}
\label{thm:martingale}Under Assumption (\ref{eq38}), we have that%
\[
\left(  \tilde{Q}^{n}(t),\tilde{X}^{n}(t),\tilde{Y}^{0,n}(t),\tilde{Y}%
^{K,n}(t)\right)  \Rightarrow\left(  \tilde{Q}(t),\tilde{X}(t),\tilde{Y}%
^{0}(t),\tilde{Y}^{K}(t)\right)  \text{, \ }as\ n\rightarrow\infty,
\]
or, in component form,%
\[
\left.
\begin{array}
[c]{c}%
\left(  \tilde{Q}_{j}^{n}(t),\tilde{Q}_{j\rightarrow
i}^{(d),n}(t),\tilde
{X}_{j}^{n}(t),\tilde{X}_{j\rightarrow i}^{(d),n}(t),\tilde{Y}_{j}%
^{0,n}(t),\tilde{Y}_{j\rightarrow i}^{0,(d),n}(t),\tilde{Y}_{j}^{K,n}%
(t)\right)  \\
\text{ \ \ \ \ \ \ \ }\Rightarrow\left(  \tilde{Q}%
_{j}(t),\tilde{Q}_{j\rightarrow i}^{(d)}(t),\tilde{X}_{j}(t),\tilde
{X}_{j\rightarrow
i}^{(d)}(t),\tilde{Y}_{j}^{0}(t),\tilde{Y}_{j\rightarrow
i}^{0,(d)}(t),\tilde{Y}_{j}^{K}(t)\right),
\end{array}
\right.  as\text{ }n\rightarrow\infty,
\]
where $\tilde{X}(t)$ is a Brownian motion with covariance matrix
$\Gamma$. Moreover, $\tilde{X}(t)-\theta t$ is a martingale with
respect to the filtration
$\mathcal{F}_{t}=\sigma(\tilde{Q}(s),\tilde{Y}^{0}(s),\tilde{Y}^{K}(s),s\leq
t)$.
\end{The}
\textbf{Proof:}. First, we define%
\begin{equation}
\tau_{+}^{n}(t)=\min\{T_{k}^{n}:T_{k}^{n}>t\}\text{ and }
\tau_{-}^{n}(t)=\max\{T_{k}^{n}:T_{k}^{n}\leq t\}.\label{eq55}%
\end{equation}%
For the station node $j\in\text{SN}$, when $\tau_{+}^{j,n}(nt)$
approximates $nt$ from its right side, we have%
\begin{align}
&  \lim_{n\rightarrow\infty}E\left[  \left|  \frac{1}{\sqrt{n}}(S_{j}^{n}%
(\tau_{+}^{j,n}(nt))-b_{j}^{n}\tau_{+}^{j,n}(nt))-\tilde{S}_{j}^{n}(t)\right|
\right]  \nonumber\\
&  =\lim_{n\rightarrow\infty}E\left[  \left|  \frac{1}{\sqrt{n}}(1-b_{j}%
^{n}(\tau_{+}^{jn}(nt))-nt)\right|  \right]  \nonumber\\
&  \leq\frac{1}{\sqrt{n}}\lim_{n\rightarrow\infty}b_{j}^{n}E\left[
\tau
_{+}^{j,n}(nt)-\tau_{-}^{j,n}(nt)\right]  \nonumber\\
&  =\lim_{n\rightarrow\infty}\frac{1}{\sqrt{n}}b_{j}^{n}E\left[  u_{j}%
^{n}(1)\right]  =0.\label{eq57}%
\end{align}
Similarly, when $\tau_{-}^{j,n}(nt)$ approximates $nt$ from its left
side, we have
\begin{equation}
\lim_{n\rightarrow\infty}E[|\frac{1}{\sqrt{n}}(S_{j}^{n}%
(\tau_{-}^{j,n}(nt))-b_{j}^{n}\tau_{-}^{j,n}(nt))-\tilde{S}_{j}^{n}(t)|]  =0.\label{eq58}%
\end{equation}
Moreover, we obtain
\begin{equation}
E[\tilde{S}_{j}^{n}(T_{k+1}^{j,n})-\tilde{S}_{j}^{n}(T_{k}^{j,n}%
)|\varsigma_{T_{k}^{j,n)}}^{n}]=\frac{1}{\sqrt{n}}\{1-b_{j}^{n}E[u_{j}%
^{n}(k+1)|\varsigma_{T_{k}^{j,n)}}^{n}]\}=0,\label{eq59}%
\end{equation}
where the filtration $\left\{  \varsigma_{t}^{n}\right\}  $ is
defined in (\ref{eq54}). Notice that for any $\left\{
\varsigma_{t}^{n}\right\} $-stopping time $T$ and any random
variable $X$ with $E\left[  \left|  X\right|
\right]  <\infty$,%
\begin{equation}
E\left[  E\left[  X\left|  \varsigma_{t}^{n}\right.  \right]  \left|
\varsigma_{t}^{n}\right.  \right]  I_{\{T>t\}}=E\left[  X\left|
\varsigma _{t}^{n}\right.  \right]  I_{\{T>t\}}=E\left[
XI_{\{T>t\}}\left|
\varsigma_{t}^{n}\right.  \right]  .\label{eq60}%
\end{equation}
Also, for each $j\in\text{SN}$ and all $s,t\geq0$,
\begin{align*}
&  E\left[  \tilde{S}_{j}^{n}(t+s)-\tilde{S}_{j}^{n}(t)\left|
\varsigma
_{t}^{n}\right.  \right]  \\
&  =E\left[  \tilde{S}_{j}^{n}(t+s)-\frac{1}{\sqrt{n}}\left(  S_{j}^{n}%
(\tau_{-}^{j,n}(n(t+s)))-b_{j}^{n}\tau_{-}^{j,n}(n(t+s))\right)
\left|
\varsigma_{t}^{n}\right.  \right]  \\
& \text{ \ \ } +E\left[  \frac{1}{\sqrt{n}}\left(  S_{j}^{n}(\tau_{+}^{j,n}(nt))-b_{j}%
^{n}\tau_{+}^{j,n}(nt)\right)  -\tilde{S}_{j}^{n}(t)\left|  \varsigma_{t}%
^{n}\right.  \right]  \\
& \text{ \ \ } -\sum_{k}E\left[  E\left[  \tilde{S}_{j}^{n}(T_{k+1}^{j,n})-\tilde{S}%
_{j}^{n}(T_{k}^{j,n})\left|  \varsigma_{T_{k}^{j,n}}^{n}\right.
\right]
I_{\left\{  nt<T_{k}^{j,n}\leq n(t+s)\right\}  }\left|  \varsigma_{t}%
^{n}\right.  \right]  .
\end{align*}
Hence, it follows from (\ref{eq57}) to (\ref{eq60}) that%
\begin{equation}
\lim_{n\rightarrow\infty}E\left[  \left|  E\left[  \tilde{S}_{j}%
^{n}(t+s)-\tilde{S}_{j}^{n}(t)\left|  \varsigma_{t}^{n}\right.
\right]
\right|  \right]  =0.\label{eq61}%
\end{equation}
For road node $j\rightarrow i$ $(i,j=1,\ldots,N\text{ with }i\neq
j)$. When we approximate $nt$ from both sides, a similar analysis to
the proof of (\ref{eq61}) for station node $j$. For all $s,t\geq0$,
we have
\begin{equation}
\lim_{n\rightarrow\infty}E\left[  \left|  E\left[
\tilde{S}_{j\rightarrow
i}^{(d),n}(t+s)-\tilde{S}_{j\rightarrow i}^{(d),n}(t)\left|  \varsigma_{t}%
^{n}\right.  \right]  \right|  \right]  =0.\label{eq62}%
\end{equation}

Next, we can set up the scaling queueing processes by mean of
(\ref{eq2}) to (\ref{eq6}) for the station nodes and of the road
nodes through the scaling processes (\ref{eq51-1}) to (\ref{eq53}),
given by:
\begin{equation}
\tilde{Q}^{n}(t)=\tilde{Q}^{n}(0)+\tilde{X}^{n}(t)+R^{0,n}\tilde{Y}%
^{0,n}(t)+R^{K,n}\tilde{Y}^{K,n}(t),\label{eq63}%
\end{equation}
and for each $n$,
$(\tilde{Q}^{n}(t),\tilde{Y}^{0,n}(t),\tilde{Y}^{K,n}(t))$ has the
properties (\ref{eq24}) to (\ref{eq30}) with the state space $S^{n}$
as follow:
\begin{align*}
S^{n}  & \equiv\left\{  x=(x_{1},\ldots,x_{N^{2}})^{^{\prime}}\in R_{+}%
^{N^{2}}:x_{i}\leq\tilde{K}_{i}^{n}=\frac{K_{i}^{n}}{\sqrt{n}}\text{ for }%
i\in\text{SN,}\right.\\
& \left.  \text{ \ \ \ \ \ and
}x_{i}\leq\frac{\sum_{i=1}^{N}C_{i}^{n}}{\sqrt{n}}+1\text{ for
}i\in\text{RN}\right\}  .
\end{align*}
For station node $j=1,\ldots,N$, by using (\ref{eq5}), (\ref{eq13}),
(\ref{eq51-1}) to (\ref{eq53})
and $\tilde{X}_{j}^{n}(t)=\frac{1}{\sqrt{n}}X_{j}^{n}(nt)=\sqrt{n}\bar{X}%
_{j}^{n}(t)$, we have%
\begin{equation}
\tilde{X}_{j}^{n}(t)=\tilde{Q}_{j}^{n}(0)+\frac{1}{\sqrt{n}}\sum_{d=1}^{2}\sum_{i\neq
j}^{N}\hat {S}_{i\rightarrow j}^{(d),n}(n\bar{B}_{i\rightarrow
j}^{(d),n}(t))-\frac{1}{\sqrt{n}}\hat{S}_{j}^{n}(n\bar{B}_{j}^{n}(t))+\sqrt{n}\theta
_{j}^{n}t.\label{eq64}
\end{equation}
For road node $j\rightarrow i\text{ }(i,j=1,\ldots,N\text{ with
}i\neq j)$, by using (\ref{eq13}), (\ref{eq17}), (\ref{eq51-1}) to
(\ref{eq53}) and $\tilde{X}_{j\rightarrow
i}^{(d),n}(t)=\frac{1}{\sqrt{n} }X_{j\rightarrow
i}^{(d),n}(nt)=\sqrt{n}\bar{X}_{j\rightarrow
i}^{(d),n}(t)$, we have,%
\begin{align}
\tilde{X}_{j\rightarrow i}^{(1),n}(t) &
=\tilde{Q}_{j\rightarrow i}^{(1),n}(0)+\frac{1}{\sqrt{n}%
}\hat{R}_{i}^{j,n}(n\bar{S}_{j}^{n}(\bar{B}_{j}(t)))\nonumber\\
& \text{ \ \ } +\frac{1}{\sqrt{n}}p_{j\rightarrow i}\hat{S}_{j}^{n}(n\bar{B}_{j}%
^{n}(t))-\frac{1}{\sqrt{n}}\hat{S}_{j\rightarrow i}^{(1),n}(n\bar
{B}_{j\rightarrow
i}^{(1),n}(t))+\frac{1}{\sqrt{n}}\theta_{j\rightarrow
i}^{(1),n}nt\label{eq65}%
\end{align}
and
\begin{align}
\tilde{X}_{j\rightarrow i}^{(2),n}(t)
=\tilde{Q}_{j\rightarrow i}^{(2),n}(0)+\frac{1}{\sqrt{n}%
}\hat{\bar{R}}_{i}^{j,n}(n\bar{Y}_{j}^{K,n}(t))-\frac{1}{\sqrt{n}}\hat{S}_{j\rightarrow
i}^{(2),n}(n\bar{B}_{j\rightarrow
i}^{(2),n}(t))+\frac{1}{\sqrt{n}}\theta_{j\rightarrow i}^{(2),n}
nt.\label{eq66}%
\end{align}
From Assumption (\ref{eq38}), using the Continuous Mapping Theorem
and Theorem \ref{thm:fluid} (Fluid Limit), we obtain that for
station node $j$,
\begin{equation}
\tilde{X}_{j}^{n}(t)\Rightarrow\tilde{X}_{j}(t)=\tilde{Q}_{j}(0)+\sum_{d=1}^{2}\sum_{i\neq
j}^{N}\tilde{S}_{i\rightarrow j}^{(d)}(t)-\tilde{S}%
_{j}(t)+\theta_{j}t,\label{eq67}%
\end{equation}
where $\tilde{X}_{j}(t)$ is an Brownian motion with the initial
queue length $\tilde{Q}_{j}(0)$ and the drift $\theta_{j}$.
For road station $j\rightarrow i$,%
\begin{equation}
\tilde{X}_{j\rightarrow
i}^{(1),n}(t)\Rightarrow\tilde{X}_{j\rightarrow
i}^{(1)}(t)=\tilde{Q}_{j\rightarrow i}^{(1)}(0)+\tilde{R}_{i}^{j}%
(b_{j}t)+p_{j\rightarrow i}\tilde{S}_{j}(t)-\tilde{S}_{j\rightarrow i}^{(1)}%
(t)+\theta_{j\rightarrow i}^{(1)}t,\label{eq69}%
\end{equation}
where $\tilde{X}_{j\rightarrow i}^{(1)}(t)$ is an Brownian motion
with the initial queue length $\tilde{Q}_{j\rightarrow i}^{(1)}(0)$
and the drift $\theta_{j\rightarrow i}^{(1)}$. Similarly we have
\begin{equation}
\tilde{X}_{j\rightarrow
i}^{(2),n}(t)\Rightarrow\tilde{X}_{j\rightarrow
i}^{(2)}(t)=\tilde{Q}_{j\rightarrow
i}^{(2)}(0)+\tilde{\bar{R}}_{i}^{j,n}(\bar{Y}_{j}^{K,n}(t))-\tilde{S}_{j\rightarrow
i}^{(2)}(t)+\theta
_{j\rightarrow i}^{(2)}t,\label{eq71}%
\end{equation}
where $\tilde{X}_{j\rightarrow i}^{(2)}(t)$ is an Brownian motion
with the initial queue length $\tilde{Q}_{j\rightarrow i}^{(2)}(0)$
and the drift $\theta_{j\rightarrow i}^{(2)}$. The covariance matrix
$\Gamma=(\Gamma_{\tilde{k},\tilde{l}} )_{N^{2}\times N^{2}}$ of
$\tilde{X}(t)=(\tilde{X}_{j}(t),$ $\tilde{X}_{j\rightarrow
i}^{(d)}(t))$ is given by%
\begin{equation}
\Gamma_{\tilde{k},\tilde{l}}=\left\{
\begin{array}
[c]{l}%
\left.
\begin{array}
[c]{l}%
\sum_{d=1}^{2}\sum_{i\neq k}^{N}b_{i\rightarrow
k}^{(d)}(c_{s,i\rightarrow
k}^{(d)})^{2}\delta_{\tilde{k},\tilde{l}}\\
\text{ \ }+b_{l}c_{a,l}^{2}\delta_{\tilde{k},\tilde{l}},
\end{array}
\right.  \text{ \ \ \ \ if }\sigma(S_{k})=\tilde{k},\sigma(S_{l})=\tilde{l};\\
\text{ \ }p_{k\rightarrow l}b_{k}c_{a,k}^{2},\text{
\ \ \ \ \ \ \ \ \ \ \ \ \ \ \ \ \ \ \  \ \ \ \ \ \ \ \ \ \ if }\sigma(S_{k}%
)=\tilde{k},\sigma(R_{k\rightarrow l})=\tilde{l},d=1;\\
\left.
\begin{array}
[c]{l}%
b_{i}p_{i\rightarrow k}(\delta_{\tilde{k},\tilde{l}}-p_{i\rightarrow l})\\
\text{ \ }+p_{k\rightarrow l}b_{k}c_{a,k}^{2}+b_{k\rightarrow
l}^{(d)}(c_{s,k\rightarrow l}^{(d)})^{2},
\end{array}
\right.  \text{ \ if }\sigma(R_{k\rightarrow
l})=\tilde{k},\tilde{l}=\tilde
{k},d=1;\\
\left.
\begin{array}
[c]{l}%
b_{k}\alpha_{i\rightarrow
k}(\delta_{\tilde{k},\tilde{l}}-\alpha_{i\rightarrow
l})\\
\text{ \ }+b_{k\rightarrow l}^{(d)}(c_{s,k\rightarrow l}^{(d)})^{2},
\end{array}
\right.  \text{ \ \ \ \ \ \ \ \ \ \ \ \ \ \ \ \ if
}\sigma(R_{k\rightarrow l})=\tilde{k},\tilde{l}=\tilde
{k},d=2;\\
\text{ \ }0,\text{
\ \ \ \ \ \ \ \ \ \ \ \ \ \ \ \ \ \ \ \ \ \ \ \ \ \ \ \ \ \ \ \ \ \ \ \ \ \ \ \ \ otherwise}%
.
\end{array}
\right.  \label{eq72}%
\end{equation}
Now, let $h(t)$ be an arbitrary real, bounded and continuous
function. For an arbitrary positive integer $m$, let $t_{i}\leq
t\leq t+s,i\leq m$. Define
\[
\tilde{H}^{n}(t)=\left(  \tilde{Q}^{n}(t),\tilde{Y}^{0,n}(t),\tilde{Y}%
^{K,n}(t)\right),  \text{ \  }\tilde{H}(t)=\left(  \tilde{Q}(t),\tilde{Y}%
^{0}(t),\tilde{Y}^{K}(t)\right)  ,
\]%
\[
G^{n}(t,s)=\left(  G_{j}^{n}(t,s),G_{j\rightarrow
i}^{(d),n}(t,s)\right)  ,
\]%
\[
G_{j}^{n}(t,s)=\tilde{X}_{j}^{n}(t+s)-\tilde{X}_{j}^{n}(t)\text{ ,
}G_{j\rightarrow i}^{(d),n}(t,s)=\tilde{X}_{j\rightarrow i}^{(d),n}%
(t+s)-\tilde{X}_{j\rightarrow i}^{(d),n}(t).
\]
Notice that%
\[
\tilde{S}_{j}^{n}(t)=\frac{1}{\sqrt{n}}\left(  \sup\left\{  k:\sum_{l=1}%
^{k}u_{j}^{n}(l)\leq b_{j}^{n}nt\right\}  -b_{j}^{n}nt\right)  ,
\]%
\[
\tilde{S}_{j\rightarrow i}^{(d),n}(t)=\frac{1}{\sqrt{n}}\left(
\sup\left\{
k:\sum_{l=1}^{k}v_{j\rightarrow i}^{(d),n}(l)\leq b_{j\rightarrow i}%
^{(d),n}nt\right\}  -b_{j\rightarrow i}^{(d),n}nt\right),
\]
by using the Assumption (\ref{eq38}), there exist some nonnegative
constants $C_{1}$ and $C_{2}$ such that $b_{j}^{n}\leq C_{1}$ and
$b_{j\rightarrow i}^{(d),n}\leq C_{2}$. From the convergences of
(\ref{eq61}) and (\ref{eq62}), we have%
\begin{align*}
&  \left|  E\left[  h\left(  \tilde{H}(t_{i}),i\leq m\right)  \left(
\tilde{X}(t+s)-\tilde{X}(t)-\theta s\right)  \right]  \right|  \\
&  =\left|  \lim_{n\rightarrow\infty}E\left[  h\left(  \tilde{H}^{n}%
(t_{i}),i\leq m\right)  G^{n}(t,s)\right]  \right|  \\
&  =\lim_{n\rightarrow\infty}\left|  E\left[  h\left(  \tilde{H}^{n}%
(t_{i}),i\leq m\right)  E\left[  G^{n}(t,s)\left|
\varsigma_{t}^{n}\right.
\right]  \right]  \right|  \\
&  \leq M\lim_{n\rightarrow\infty}E\left[  \left|  E\left[
G^{n}(t,s)\left|
\varsigma_{t}^{n}\right.  \right]  \right|  \right]  \\
&  =0,
\end{align*}
where $M$ is some positive constant. The arbitrariness of $h(t)$,
$t_{i}$, $t$
and $t+s$ implies that%
\[
E\left[  \tilde{X}(t+s)-\tilde{X}(t)-\theta s\left|
\mathcal{F}_{u},u\leq t\right.  \right]  =0.
\]
This shows that $\tilde{X}(t)-\theta t$ is an $\left\{  \mathcal{F}%
_{t}\right\}  $-martingale. This completes the proof.
\textbf{{\rule{0.08in}{0.08in}}}

\begin{Rem}
\label{rem:citeDaiW}Note that Dai \cite{DaiW:1996} discussed the
queueing networks with finite buffers, this paper is related well to
fluid and diffusion limits in Dai \cite{DaiW:1996} in order to deal
with a two-class closed queueing network.
\end{Rem}

Now, we give the diffusion limit for the bike sharing system. In
Section \ref{sec5:fluid limit}, we set up a sequence of closed
queueing networks corresponding to the bike sharing systems, and
prove the limit theorems of the fluid scaled equations of the busy
period processes and the idle period processes through the
functional strong law of large numbers and the oscillation property
of an $(S^{n},R^{n})$-regulation. This is summarized as the Fluid
Limit Theorem \ref{thm:fluid}. Furthermore, based on the Fluid Limit
Theorem, we prove the weak limit of the diffusion scaled processes
of some performance measures and obtain a key martingale. Also see
Theorem \ref{thm:martingale}.

The following theorem provides a diffusion limit, and its proof is
easy by means of some similar analysis to Theorems 3.2 and 3.3 in
Dai \cite{DaiW:1996} or Theorem 3.1 in Dai and Dai \cite{Dai:1999}.

\begin{The}
\label{thm:diffusion}(Diffusion Limit Theorem) Under
Assumption (\ref{eq38}), we have%
\[
\left(  \frac{1}{\sqrt{n}}Q^{n}(nt),\frac{1}{\sqrt{n}}Y^{0,n}(nt),\frac{1}%
{\sqrt{n}}Y^{K,n}(nt)\right)  \Rightarrow\left(  \tilde{Q}(t),\tilde{Y}%
^{0}(t),\tilde{Y}^{K}(t)\right),
\]
where $\tilde{Q}(t)=\left(  \tilde{Q}_{j}(t),\tilde{Q}_{j\rightarrow i}%
^{(d)}(t)\right)  $, $\tilde{Y}^{0}(t)=\left(
\tilde{Y}_{j}^{0}(t),\tilde {Y}_{j\rightarrow i}^{0,(d)}(t)\right)
$; $\tilde{Q}(t)$ together with ${\tilde{Y}}^{0}{(t)}$ and
$\tilde{Y}^{K}(t)$\ are an $(S,\theta,\Gamma,R)$-semimartingale
reflecting Brownian motion
with $\tilde{Q}(t)=\tilde{Q}(0)+\tilde{X}(t)+R^{0}\tilde{Y}%
^{0}(t)+R^{K}\tilde{Y}^{K}(t)$. The state space $S$ is given by
(\ref{eq31}) to (\ref{eq32-2}). For station node $j$,
$\tilde{X}_{j}(t)$ is given by (\ref{eq67}), $R^{0}$ and $R^{K}$ are
given by (\ref{eq11}), (\ref{eq12}). For road node $j\rightarrow i$,
when $d=1$, $\tilde{X}_{j\rightarrow i}^{(1)}(t)$ is given by
(\ref{eq69}), $R^{0}$ and $R^{K}$ are given by (\ref{eq15}) and
(\ref{eq16}); when $d=2$, $\tilde{X}_{j\rightarrow i}^{(2),n}(t)$ is
given by (\ref{eq71}), $R^{0}$ and $R^{K}$ are given by (\ref{eq19})
and (\ref{eq20}), and the covariance matrix
$\Gamma=(\Gamma_{\tilde{k},\tilde{l}} )_{N^{2}\times N^{2}}$ of
$\tilde{X}(t)=(\tilde{X}_{j}(t),\tilde{X}_{j\rightarrow
i}^{(d)}(t))$ is given by (\ref{eq72}).
\end{The}

\section{Performance analysis}\label{sec7:performance}

In this section, we first set up a basic adjoint relationship for
the steady-state probabilities of $N$ station nodes and of $N(N-1)$
road nodes in the multiclass closed queueing network. Then we
analyze some key performance measures of the bike sharing system.

From Theorem \ref{thm:diffusion}, it is seen that the scaling
queueing processes, for the numbers of bikes in the stations and on
the roads, converge in distribution to a semimartingale reflecting
Brownian motion
$\tilde{Q}(t)=\left({\tilde{Q}_{\tilde{i}}(t),\tilde{Q}_{\tilde{j}}^{(d)}(t)}\right)$
for $\tilde{i}=\sigma(S_{i})$ $(i=1,\ldots,N)$ and
$\tilde{j}=\sigma(R_{j\rightarrow i})$ $(i,j=1,\ldots,N\text{ with
}i\neq j,d=1,2)$, where the state space $S$, the drift vector
$\theta=\left({\theta_{\tilde{i}},\theta_{\tilde{j}}^{(d)}}\right)$
for $\tilde{i}=\sigma({S_{i}}),\tilde{j}=\sigma({R_{j\rightarrow
l}})$, the covariance matrix
\begingroup
\renewcommand*{\arraystretch}{1.5}
\[
\Gamma=\left(
\begin{array}
[c]{cc}%
(  \Gamma_{\tilde{i},\tilde{k}})   & (\Gamma
_{\tilde{i},\tilde{j}}^{(d)})  \\
(\Gamma_{\tilde{l},\tilde{k}}^{(d)})   & (\Gamma
_{\tilde{l},\tilde{j}}^{(d)})
\end{array}
\right)  _{N^{2}\times N^{2}}%
\]
\endgroup
for
$\tilde{i}=\sigma(S_{i}),\tilde{k}=\sigma(S_{k}),\tilde{j}=\sigma(R_{j\rightarrow
h}),\tilde{l}=\sigma(R_{l\rightarrow g})$ and the reflecting matrix
$R=\left({\left({R_{\tilde{i}}^{0}}\right),\left({R_{\tilde{j}}^{K,(d)}}\right)}\right)$
for $\tilde{i}=\sigma({S_{i}}),\tilde{j}=\sigma({R_{j\rightarrow
l}})$, as seen in those previous sections. Hence, it is natural to
approximate the steady-state distribution of the queue-length
process by means of the steady-state distribution of the
semimartingale reflecting Brownian motion.

From Lemma \ref{lem:promatrix} and Theorem 1.3 in Dai and Williams
\cite{Dai:1996}, it is seen that there exists a unique stationary
distribution
$\pi=\left({\pi_{\tilde{i}},\pi_{\tilde{j}}^{(d)}}\right)$ on
$(S,\mathcal {B}_{S})$ for the SRBM
$\tilde{Q}(t)=\left({\tilde{Q}_{\tilde{i}}(t),\tilde{Q}_{\tilde{j}}^{(d)}(t)}\right)$.
Furthermore,
$\pi=\left({\pi_{\tilde{i}},\pi_{\tilde{j}}^{(d)}}\right)$ is
equivalent to the Lebesgue measure on the state space $S$, thus for
every bounded Borel function $f$ on $S$ and for $t\geq0$, we have
\[
E_{\pi}\left[f\left({\tilde{Q}(t)}\right)\right]\equiv\int_{S}\left(
{E_{x}\left[f\left({\tilde{Q}(t)}\right)\right]}\right)\pi(dx)=\int_{S}f(x)\pi(dx).
\]
Then for each $\tilde{i}=1,\ldots,N$ (i.e.,
$\tilde{i}=\sigma(S_{i}),i=1,\ldots,N$) and
$\tilde{j}=1,\ldots,N(N-1)$ (i.e., $\tilde{j}=\sigma(R_{j\rightarrow
i}),i,j=1,\ldots,N\text{ with }i\neq j$), let
$\delta=\left({\delta_{\tilde{i}},\delta_{\tilde{j}}^{(d)}}\right)$
denote $(N^{2}-1)$-dimensional Lebesgue measure (surface measure)
vector on face $(F,\mathcal{B}_{F})$. Thus, there is a finite Borel
measure vector
$\beta^{F}=\left({\beta_{\tilde{i}}^{F},\beta_{\tilde{j}}^{F,(d)}}\right)$
on face $F=(F_{\tilde{i}},F_{\tilde{j}})$ such that
$\beta^{F}\approx\delta$ and
\[
E_{\pi}\left\{\int_{0}^{t}1_{A}\left({\tilde{Q}(s)}\right)d\tilde{Y}(s)\right\}=t\beta^{F}(A),\text{
\ }t\geq0,A\in\mathcal{B}_{F},
\]
where
$\tilde{Y}(t)=\left({\tilde{Y}^{0}(t),\tilde{Y}^{K}(t)}\right)$.
Notice that the SRBM
$\tilde{Q}(t)=(\tilde{Q}_{\tilde{i}}(t),\tilde{Q}_{\tilde{j}}^{(d)}(t))$
is a strong Markov process with continuous sample paths.
Furthermore, let
$p(x)=\left({p_{\tilde{i}}(x_{\tilde{i}}),p_{\tilde{j}}^{(d)}\left({x_{\tilde{j}}^{(d)}}\right)}\right)$,
$p^{F}(x)=\left({p_{\tilde{i}}^{F}\left({\delta_{\tilde{i}}}\right),p_{\tilde{j}}^{F,(d)}\left({\delta_{\tilde{j}}^{(d)}}\right)}\right)$,
and define $d\pi=pdx$, i.e.,
$d\pi_{\tilde{i}}=p_{\tilde{i}}dx_{\tilde{i}}$ for
$\tilde{i}=\sigma(S_{i})$ $(i=1,\ldots,N)$ and
$d\pi_{\tilde{j}}^{(d)}=p_{\tilde{j}}^{(d)}dx_{\tilde{j}}^{(d)}$ for
$\tilde{j}=\sigma(R_{j\rightarrow i})$ $(i,j=1,\ldots,N\text{ with
}i\neq j,d=1,2)$. Further, we define $d\beta^{F}=p^{F}d\delta$,
i.e., $d\beta_{\tilde{i}}^{F}=p_{\tilde{i}}^{F}d\delta_{\tilde{i}}$
for $\tilde{i}=\sigma(S_{i})$ $(i=1,\ldots,N)$ and
$d\beta_{\tilde{j}}^{F,(d)}=p_{\tilde{j}}^{F,(d)}d\delta_{\tilde{j}}^{(d)}$
for $\tilde{j}=\sigma(R_{j\rightarrow i})$ $(i,j=1,\ldots,N\text{
with }i\neq j,d=1,2)$. Let $\nabla f(x)$ be the gradient of $f$, and
$C_{b}^{2}(S)$ the space of twice differentiable functions whose
first and second order partial derivative are continuous and bounded
on the state space $S$. Base on this, it follows from the Ito's
formula that the probability measures $p(x)$ and $p^{F}(x)$ have a
basic adjoint relationship as follows: for $\forall f\in
C_{b}^{2}(S)$,
\begin{align}
&  \int_{S}\left(  {\mathcal{L}}f{(x)}p(x)\right)
dx+\sum_{\tilde{i}=1}^{N}\int_{F_{\tilde{i}}^{\blacktriangledown}}(\mathcal{D}_{\tilde{i}}f(\delta_{\tilde
{i}})p_{\tilde{i}}^{F}(\delta_{\tilde{i}}))d\delta_{\tilde{i}}+\sum_{d=1}^{2}\sum_{\tilde{j}=1}^{N^{2}-N}\int_{F_{\tilde{j}}^{\blacktriangledown}
}(\mathcal{D}_{\tilde{j}}f(\delta_{\tilde{j}}^{(d)})p_{\tilde{j}}^{F,(d)}%
(\delta_{\tilde{j}}^{(d)}))d\delta_{\tilde{j}}^{(d)}\nonumber\\
&  \text{ \
}+\sum_{\tilde{i}=1}^{N}\int_{F_{\tilde{i}}^{\blacktriangle}}(\mathcal{D}_{\tilde{i}}f(\delta_{\tilde
{i}})p_{\tilde{i}}^{F}(\delta_{\tilde{i}}))d\delta_{\tilde{i}}+\sum_{d=1}^{2}\sum_{\tilde{j}=1}^{N^{2}-N}\int_{F_{\tilde{j}}^{\blacktriangle}
}(\mathcal{D}_{\tilde{j}}f(\delta_{\tilde{j}}^{(d)})p_{\tilde{j}}^{F,(d)}%
(\delta_{\tilde{j}}^{(d)}))d\delta_{\tilde{j}}^{(d)}=0,\label{PA-eq1}%
\end{align}
where
\[
\mathcal{L}f =
\sum_{\tilde{i}=1}^{N}\mathcal{L}f(x_{\tilde{i}})+\sum_{d=1}^{2}\sum_{\tilde{j}}^{N^{2}-N}\mathcal{L}f(x_{\tilde{j}}^{(d)}),
\]
for $i,k,j=1,\ldots, N\text{ with }i\neq j,d=1,2$, and
$\tilde{k}=\sigma(S_{k}),\tilde{j}=\sigma(R_{j\rightarrow i})$,
$\tilde{i}=\sigma(S_{i})\in\{1,2,\ldots,N\}$,

\[
\mathcal{L}f(x_{\tilde{i}})=\frac{1}{2}\sum_{\tilde{k}=1}^{N}\Gamma
_{\tilde{i},\tilde{k}}\frac{\partial^{2}f(x_{\tilde{i}})}{\partial
x_{\tilde{i}}\partial
x_{\tilde{k}}}+\frac{1}{2}\sum_{d=1}^{2}\sum_{\tilde{j}=1}^{N^{2}-N}\Gamma
_{\tilde{i},\tilde{j}}^{(d)}\frac{\partial^{2}f(x_{\tilde{i}})}{\partial
x_{\tilde{i}}\partial
x_{\tilde{j}}^{(d)}}+\theta_{\tilde{i}}\frac{\partial
f(x_{\tilde{i}})}{\partial x_{\tilde{i}}},
\]%
\[
\mathcal{D}_{\tilde{i}}f(\delta_{\tilde{i}})\equiv
v_{\tilde{i}}^{\prime}\nabla
f(\delta_{\tilde{i}})=\sum_{\tilde{k}=1}^{N}v_{\tilde{k},\tilde{i}}\frac{\partial}{\partial
\delta_{\tilde{k}}}f(\delta_{\tilde{i}})+\sum_{d=1}^{2}\sum_{\tilde{j}=1}^{N^{2}-N}v_{\tilde{j},\tilde{i}}\frac{\partial}{\partial
\delta_{\tilde{j}}^{(d)}}f(\delta_{\tilde{i}}),
\]
for $l,k,i,j,h=1,\ldots, N\text{ with }j\neq i,l\neq k\text{ and
}d=1,2$, and $\tilde{l}=\sigma(R_{l\rightarrow
k}),\tilde{h}=\sigma(S_{h})$, $\tilde{j}=\sigma(R_{j\rightarrow
i})\in\{1,2,\ldots,N^{2}-N\}$,

\[
\mathcal{L}f(x_{\tilde{j}}^{(d)})=\frac{1}{2}\sum_{\tilde{h}=1}^{N}\Gamma
_{\tilde{j},\tilde{h}}^{(d)}\frac{\partial^{2}f(x_{\tilde{j}}^{(d)})}{\partial
x_{\tilde{j}}^{(d)}\partial
x_{\tilde{h}}}+\sum_{\tilde{l}=1}^{N^{2}-N}\Gamma
_{\tilde{j},\tilde{l}}^{(d)}\frac{\partial^{2}f(x_{\tilde{j}}^{(d)})}{\partial
x_{\tilde{j}}^{(d)}\partial x_{\tilde{l}}^{(d)}}
+\theta_{\tilde{j}}^{(d)}\frac{\partial
f(x_{\tilde{j}}^{(d)})}{\partial x_{\tilde{j}}^{(d)}},
\]

\[
\mathcal{D}_{\tilde{j}}f(\delta_{\tilde{j}}^{(d)})\equiv
v_{\tilde{j}}^{\prime}\nabla
f(\delta_{\tilde{j}}^{(d)})=\sum_{\tilde{l}=1}^{N(N-1)}v_{\tilde{l},\tilde{j}}\frac{\partial}{\partial
\delta_{\tilde{l}}^{(d)}}f(\delta_{\tilde{j}}^{(d)})+\sum_{\tilde{h}=1}^{N}v_{\tilde{h},\tilde{j}}\frac{\partial}{\partial
\delta_{\tilde{h}}}f(\delta_{\tilde{j}}^{(d)}),
\]
$F_{\tilde{i}}^{\blacktriangledown}$ and
$F_{\tilde{i}}^{\blacktriangle}$ denote the ``bottom face'' and the
``top face'' in this box state space $S$ corresponding to empty
station $i$ and full station $i$, respectively. As a similar
expression, it is clear that $F_{\tilde{j}}^{\blacktriangledown}$
and $F_{\tilde{j}}^{\blacktriangle}$ are related to road
$j\rightarrow i$; $v_{\tilde{k}}$ is the $\tilde{k}$th column of the
reflection matrix
$R=\left({\left({R_{\tilde{i}}^{0}}\right),\left({R_{\tilde{j}}^{K,(d)}}\right)}\right)$.

Now, we consider some key performance measures of the bike sharing
system in terms of the steady-state probability density function $p$
on $(S,\mathcal {B}_{S})$ and an nonnegative integrable Borel
function $p^{F}$ on $(F,\mathcal {B}_{F})$. Here, it is easy to see
that for $\tilde{i}=1,\ldots,N$ and $\tilde{j}=1,\ldots,N(N-1)$, the
``bottom face'' $F_{\tilde{i}}^{\blacktriangledown}$
($F_{\tilde{j}}^{\blacktriangledown}$) and the ``top face''
$F_{\tilde{i}}^{\blacktriangle}$ ($F_{\tilde{j}}^{\blacktriangle}$)
are precisely parallel in this box state space $S$.

{\bf(1)} The steady-state probability that station $i$ is empty is
given by
\[
\int_{S}p_{\tilde{i}}^{F}1_{\{x_{\tilde{i}}\in
F_{\tilde{i}}^{\blacktriangledown}\}}dx_{\tilde{i}}, \text{\ \ \ \ \
\ \ \ for } \tilde{i}=\sigma(S_{i}).
\]

{\bf(2)} The steady-state probability that station $i$ is full is
given by
\[
\int_{S}p_{\tilde{i}}^{F}1_{\{x_{\tilde{i}}\in
F_{\tilde{i}}^{\blacktriangle}\}}dx_{\tilde{i}}, \text{\ \ \ \ \ \ \
\ for } \tilde{i}=\sigma(S_{i}).
\]

{\bf(3)} The steady-state probability that road $j\rightarrow i$ is
empty for bikes of class $d$ is given by
\[
\int_{S}p_{\tilde{j}}^{F,(d)}1_{\{x_{\tilde{j}}^{(d)}\in
F_{\tilde{j}}^{\blacktriangledown}\}}dx_{\tilde{j}}^{(d)}, \text{\ \
\ \ for } \tilde{j}=\sigma(R_{j\rightarrow i}),d=1,2.
\]

{\bf(4)} The steady-state probability that road $j\rightarrow i$ is
full for bikes of class $d$ is given by
\[
\int_{S}p_{\tilde{j}}^{F,(d)}1_{\{x_{\tilde{j}}^{(d)}\in
F_{\tilde{j}}^{\blacktriangle}\}}dx_{\tilde{j}}^{(d)}, \text{\ \ \ \
for } \tilde{j}=\sigma(R_{j\rightarrow i}),d=1,2.
\]

{\bf(5)} The steady-state means of the number of bikes parked at the
station $i$ and the number of bikes of class $d$ ridden on road
$j\rightarrow i$ are respectively given by
\[
\mathcal{Q}_{\tilde{i}}=\int_{S}x_{\tilde{i}}p_{\tilde{i}}(x_{\tilde{i}})dx_{\tilde{i}},\text{
\ for }\tilde{i}=\sigma(S_{i}),
\]
\[
\mathcal{Q}_{\tilde{j}}^{(d)}=\int_{S}x_{\tilde{j}}^{(d)}p_{\tilde{j}}^{(d)}\left({x_{\tilde{j}}^{(d)}}\right)dx_{\tilde{j}}^{(d)},\text{
\ for }\tilde{j}=\sigma(R_{j\rightarrow i}),d=1,2.
\]

{\bf(6)} The steady-state mean of the number of bikes of class $d$
deflecting from the full station $i$ is given by
\[
\mathcal {E}_{\tilde{i}}^{(d)}=\int_{F_{\tilde{i}}^{\blacktriangle}}
x_{\tilde{i}}^{(d)}p_{\tilde{i}}^{F,(d)}\left({x_{\tilde{i}}^{(d)}}\right)dx_{\tilde{i}}^{(d)},\text{
\ for }\tilde{i}=\sigma(S_{i}),d=1,2.
\]

\section{Concluding Remarks}\label{sec8:conclusions}
In this paper, we describe a more general large-scale bike sharing
system having renewal arrival processes and general travel times,
and develop fluid and diffusion approximation of a multiclass closed
queuing network which is established from the bike sharing system
where bikes are regarded as virtual customers, and stations and
roads are viewed as virtual nodes or servers. From the multiclass
closed queuing network, we show that the scaling queue-length
processes, which are set up by means of the number of bikes both at
stations and on roads, converge in distribution to a semimartingale
reflecting Brownian motion. Also, we obtain the Fluid Limit Theorem
and the Diffusion Limit Theorem. Based on this, we provide
performance analysis of the bike sharing system. Therefore, the
results of this paper give new highlight in the study of more
general large-scale bike sharing systems. The methodology developed
here can be applicable to deal with more general bike sharing
systems by means of the fluid and diffusion approximation. Along
such a line, there are some interesting directions in our future
research, for example,

\begin{itemize}

\item analyzing bike repositioning policies through several fleets of trucks under information technologies;

\item making price regulation of bike sharing systems through Brownian approximation of multiclass closed queuing network;

\item developing heavy traffic approximation for time-varying or periodic bike sharing systems; and

\item developing heavy traffic approximation for new ride sharing (bike or car) systems with scheduling, matching and control.

\end{itemize}

\section*{Acknowledgments}
This work was supported in part by the National Natural Science
Foundation of China under grant No. 71671158 and No. 71471160, and
Natural Science Foundation of Hebei under grant No. G2017203277.

\end{document}